\RequirePackage{fix-cm}
\documentclass[smallextended]{svjour3}       
\smartqed  
\usepackage{graphicx}
\RequirePackage[numbers]{natbib}
\usepackage{multirow}
\usepackage[ruled,vlined]{algorithm2e}
\usepackage{graphicx}
\usepackage{setspace}
\usepackage{float}
\usepackage{breakcites}
\usepackage{epsfig}
\usepackage{amssymb}
\usepackage{amsmath}
\usepackage{hyperref}
\usepackage{verbatim}
\usepackage{subcaption}
\usepackage{bm}
\usepackage{tikz}
\usepackage{color}
\newcommand{\AZ}{\textcolor{black}}

\usepackage{bbm}
\usetikzlibrary{arrows}
\usepackage[margin=1.2in]{geometry}

\tikzset{cross/.style={cross out, draw=black, minimum size=2*(#1-\pgflinewidth), inner sep=0pt, outer sep=0pt},
cross/.default={5pt}}
\usetikzlibrary{shapes.misc}
\usepackage{array}
\newcommand{\be}{\begin{eqnarray}}
\newcommand{\ee}{\end{eqnarray}}
\newcommand{\bea}{\begin{eqnarray*}}
\newcommand{\eea}{\end{eqnarray*}}

\usepackage{nicefrac}
\numberwithin{equation}{section}

 \newcolumntype{P}[1]{>{\centering\arraybackslash}m{#1}}

\newcommand{\X}{{\cal X}}

\begin{document}

\title{
Improving exploration strategies in large dimensions and rate of convergence of global random search algorithms
}

\titlerunning{Exploration strategies in large dimensions}


\author{  Jack Noonan and {Anatoly Zhigljavsky~(Corresponding author)}
}


\institute{
}

\date{}
\maketitle



\raggedbottom

\abstract{
We consider global optimization problems, where the feasible region $\X$ is a compact subset of $\mathbb{R}^d$ with $d \geq 10$. For these problems, we demonstrate that   the actual convergence of global random search algorithms is much slower than that given by the classical estimates, based on the asymptotic properties of random points, and 
that  the usually recommended space exploration schemes  are  inefficient in  the non-asymptotic regime. Moreover, we show that  uniform sampling on entire~$\X$ is much less efficient than uniform sampling on a suitable subset of $\X$,  and that the effect of replacement of  random points by low-discrepancy sequences can be felt in small dimensions only.
 }

\section{Introduction}

\label{sec:intro}

Consider the general  problem of continuous global minimization   \mbox{$f(x)\!\rightarrow\!
{\textrm{min}}_{x\in \X}$} with objective function $f(\cdot)$ and
feasible region $\X$, which is assumed to be a compact subset of $\mathbb{R}^d$ with  ${\rm vol}(\X)>0$.
 In order to avoid unnecessary technical difficulties, we assume that $\X$ is convex.
In all numerical examples, we use $\X=[0,1]^d$.

Any global optimization algorithm combines two key strategies: exploration and exploitation. Performing exploration is equivalent to what we call ``space-filling"; that is,  choosing points which are well-spread in $\cal X$. Exploitation strategies use local information about $f$ (and perhaps derivatives of $f$) and differ greatly for different types of global optimization algorithms. In this paper, we are only concerned with the exploration stage.
 Although many of our finding can be generalized to other space-filling schemes (where space-filling is not random and the space-filling strategy changes in the course of receiving more information about the objective function),
 in this paper we concentrate on simple exploration schemes like pure random search, where  space-filling is performed by
  covering $\cal X$ with balls of given radius centered at the chosen points.
  Moreover, we assume that the points chosen at the exploration stage are independent.
  That is, we associate the exploration stage with a global random search (GRS) algorithm   producing a sequence of   random points $x_1,x_2,\ldots, x_n$,
where each point $x_j \in \X$
  has some
probability distribution $P_j$ (we write this $x_j \sim P_j$) and $x_1,x_2,\ldots, x_n$ are independent.
The value $n$ is determined by a stopping rule. We assume that  $1 \leq n_{\min} \leq n \leq n_{\max} < \infty$,
where $n_{\min} $ and $  n_{\max}$ are two given numbers. The number $  n_{\max}$ determines the maximum number of function evaluations at the exploration stage and the fact that
$  n_{\max}< \infty $ determines what we call ``the non-asymptotic regime''. In the numerical study of Section~\ref{sec:delta1}, we also use Sobol's sequence, the most widely used  low-discrepancy sequence.

\smallskip

We distinguish between `small', `medium' and `high' dimensional problems depending on the following relations between $d$ and $n_{\max}$:

\begin{itemize}
  \item[(S)] small dimensions: $n_{\min} \geq 2^d$, $n_{\max} \gg 2^d$ (hence, $\log_2 n_{\max} \gg d$);
  \item[(M)] medium dimensions: $n_{\max} $ is comparable to $ 2^d$:  $c_1 d \leq \log_2 n_{\max} \leq c_2 d$ with suitable constants $c_1$ and $c_2$: $0 \ll c_1 \leq  1 \leq c_2 \ll \infty$;
  \item[(H)] high dimensions: $n_{\max} \ll 2^d$.
\end{itemize}

Of course, there are in-between situations and the classification above depends on the cost of function evaluation.
In case of non-expensive observations and $10^3 \leq n_{\max} \leq 10^6 $, typical values of $ d$  in the three cases are:
$
\mbox{ (S):  $d \leq 10$; (M):  $10 \leq d \leq 20$; (H): $d>20$.}
$
%
%
Values of $d \approx 10 $ are border-line cases between (S) and (M) whereas $d \approx 20 $ are border-line cases between (M) and~(H).

In this study, we leave out the situation (S) of small dimensions and concentrate on situations (M) and (H). The reasons why we are not interested in the  situation (S) of small dimensions are: (a) there are too many exploration schemes available in literature in the case of small dimensions, and (b) we are interested in the situations when the asymptotic regime is out of reach, and these are the situations (M) and~(H).

In all considerations below we assume that the aim of the exploration stage is to reach a neighbourhood of an unknown point $x_* \in \X$ with high probability $\geq 1-\gamma$ (with some $\gamma>0$).
We assume that $x_*$ is uniformly distributed in $\X$ and by a neighbourhood of $x_*$ we mean the ball $B={\cal B}(x_*,\varepsilon)$ with suitable $\varepsilon>0$. In other words,  we will be interested in the problem of construction of weak coverings defined as follows.

Let $x_1, \ldots, x_n$ be some points in $\mathbb{R}^d$.
Denote
$X_n= \{x_1, \ldots, x_n\}$ and
\be
\label{eq:covering}
B(X_n,r)= \bigcup_{i=1}^n B(x_i,r)\, ,
\ee
where $r>0$ is the radius of the balls ${\cal B}(x_i,r)$ and $B(x_i,r) = \X \cap {\cal B}(x_i,r)$.
We will call $B(X_n,r)$ weak (or approximate) covering of $\X$ of level $1-\gamma$ if ${\rm vol}(B(X_n,r))/{\rm vol}(\X)\geq 1-\gamma$.

If $\gamma=0$ then $B(X_n,r)$ would make a full (strong) covering of $\X$. As demonstrated in \cite{us,second_paper,noonan2022efficient}, for any $n$ and any given $\gamma>0$, one can construct weak coverings of $\X$ with significantly smaller radii $r$ than for the case $\gamma=0$ (assuming that $d$ is not too small). This is the main reason why we are not interested in strong coverings. The second reason is that numerically checking whether the set \eqref{eq:covering} makes a full covering (for a generic $X_n$) is extremely hard in situations (M) and (H) whereas simple Monte-Carlo gives very accurate estimates of $\gamma$ for weak coverings, even for very high dimensions. For a short discussion concerning full covering and its role in optimization, see Section~\ref{sec:cover}.

The main technique of construction of weak coverings will be  generation of independent random points $x_1, \ldots, x_n$ in $\X$ with $x_j \sim P$, where $P$ is a distribution concentrated either on the whole $\X$ or a subset of $\X$.  It follows from Proposition 3.2.3 in \cite{borodachov2019discrete} that using points outside $\X$ for construction of coverings is not beneficial when $\X$ is convex and hence we will always assume that $x_j \in \X$ for all $j$.

\medskip

The following are the main messages  of the paper.
\begin{enumerate}
\item
Classical results on  convergence rates of GRS algorithms  are based on the asymptotic properties of random points uniformly distributed in $\X$; see Section~\ref{sec:classical}. In the non-asymptotic regime, however, these results give estimates on the convergence rates which are far too optimistic. We show  in Section~\ref{sec:non_classical} that for medium and high dimensions, the actual  convergence rate of GRS algorithms is much slower.
\item
The usually recommended sampling schemes (these schemes are  based on the asymptotic properties of random points) are  inefficient in  the non-asymptotic regime. In particular, as shown in Section~\ref{sec:delta1}, uniform sampling on entire~$\X$ is much less efficient than uniform sampling on a suitable subset of $\X$ \AZ{(we will refer to this phenomena as the `$\delta$-effect').}
 \item  In situations (M) and (H), the effect of replacement of  random points by low-discrepancy sequences  is negligible; see Section~\ref{sec:delta1_2}.
\end{enumerate}
We also make certain practical recommendations concerning the best exploration schemes in the situations (M) and (H) in the case $\X=[0,1]^d$.
Our main  recommendations will concern the situation (M) of medium dimensions, which we consider as the hardest for analysis. The situation (H)  is simpler than (M)  in the sense that the optimization problems in case (H) are so complicated  that very simple space-filling schemes outlined in Section~\ref{sec:conclusions} provide relatively effective sampling schemes.

\medskip

The structure of the paper is as follows. In Section~\ref{sec:classical}, which  contains no new results, we discuss the importance of covering and  review classical results on convergence and rate of convergence of general GRS algorithms. The purpose Section~\ref{sec:non_classical} is to demonstrate that
for medium and high dimensions the asymptotic regime is unachievable, and hence the actual  convergence rate of GRS algorithms is much slower than the classical estimates of the rate of convergence indicate. In Section~\ref{sec:delta1} we compare several exploration strategies and show that standard recommendations (such as: ``use a low-discrepancy sequence") are inaccurate (for medium and high dimensions). In Section~\ref{Petrov_section}, we develop accurate approximations for the volume of intersection of a cube and a ball (with arbitrary centre and any radius). The approximations of Section~\ref{Petrov_section} are used throughout numerical studies of Sections~\ref{sec:non_classical} and~\ref{sec:delta1}. In Section~\ref{sec:conclusions} we summarize our findings and give recommendations on how to perform exploration of $\X$ in medium and large dimensions.

 \section{Importance of covering and classical results on convergence and rate of convergence  of  GRS algorithms}

\label{sec:classical}

\subsection{Covering radius}

\label{sec:cover}

Consider    $X_n=\{x_1, \ldots, x_n\} $,  a set of $n$ points in $\X$.
The covering radius of $\X$ for  $X_n$ is
$ \label{eq:CR}
{\rm CR} (X_n) = \max_{x\in\X}   \rho (x,X_n)  ,
$
where
\be \label{eq:CR5}
\rho (x,X_n)= \min_{x_j\in X_n} \rho (x,x_j)\,
\ee
is the distance between a point $x \in \X$ and the point set  $X_n$.
Covering radius is also
the smallest  $r \geq 0$ such that the union of the balls with centers at $x_j \in X_n$ and radius~$r$ fully covers $\X$; that is,
$ \label{eq:CR1}
{\rm CR} (X_n)= \mbox{$\min_{ {r>0} }$} \; \mbox{ such that }  \X \subseteq {\cal B} (X_n,r)\, ,
$
where ${\cal B} (X_n,r)= \bigcup_{j=1}^n {\cal B} (x_j,r)$ and
$
\label{eq:CR2} {\cal B} (x,{ r })= \{ z \in \mathbb{R}^d :\; \rho (x,z) \leq { r } \}$
is the ball of radius $r$ and centre $x\in  \mathbb{R}^d$.
Optimal $n$-point covering is the point set $X_n^*$
such that $ {\rm CR}(X_n^*)=\min_{X_n}  {\rm CR}(X_n). $ \AZ{Most of the general considerations in the paper are valid for a general distance $\rho$, but all numerical studies are conducted for the Euclidean distance only. We will thus assume that the distance $\rho$ is Euclidean. }

Other common names for the covering radius are:  fill distance (in approximation theory; see \cite{schaback2006kernel,wendland2004scattered}),  dispersion (in Quasi Monte Carlo; see \cite[Ch. 6]{Niederreiter}),
  minimax-distance criterion (in computer experiments;    see \cite{pronzato2012design,santner2003design}) and
 coverage threshold (in probability theory; see \cite{penrose2021random}).

Point sets with small covering radius are very desirable in theory and practice  of global optimization and many branches of numerical mathematics.
In particular, the celebrated results of A.G.Sukharev  imply that any $n$-point optimal covering design $X_n^*$  provides the following:
(a)
 min-max $n$-point global optimization method in the set of all adaptive $n$-point optimization strategies, see \cite{Sukh1} and \cite[Ch.4,Th.2.1]{sukharev2012minimax},
 (b) worst-case  $n$-point multi-objective global optimization method in the set of all adaptive $n$-point algorithms, see
 \cite{vzilinskas2013worst}, and
  (c) the $n$-point min-max optimal quadrature, see \cite[Ch.3,Th.1.1]{sukharev2012minimax}.
  In all three cases, the class of (objective) functions is the class of Liptshitz functions, and
    the optimality of the design is independent of the value of the Liptshitz constant. Sukharev's results on $n$-point min-max optimal quadrature formulas  have been generalized in \cite{pages1998space}
    for functional classes different from the class of Liptshitz functions; see also  formula (2.3) in  \cite{du1999centroidal}.

 \subsection{Convergence  of  a general GRS algorithm}

\label{sec:conv}

Consider the general  problem of continuous global minimization   \mbox{$f(x)\!\rightarrow\!
{\textrm{min}}_{x\in \X}$}.
Assume that $f_*=\inf_{x \in \X} f(x)>-\infty$ and $f(\cdot)$  is continuous  at all points $x \in W
(\delta)$ for some $\delta>0$, where
$W
(\delta)\!=\!\left\{x\in \X\colon f(x)\!-\!f_* \! \leqslant \!\delta\right\}$.
 That is, we assume that $f(\cdot)$  is continuous in the neighbourhood of the set $\X_*=\left \{x_* \in \X\colon f(x_* )=f_* \right \}$   of global minimizers of $f(\cdot)$, which is non-empty but may contain more than one point $x_*$. To avoid technical difficulties, we  assume that there are only a finite number of global minimizers of $f(\cdot)$; that is, the set $\X_*$ is finite.

Consider a general GRS algorithm   producing a sequence of   random points $x_1,x_2,\ldots$,
where each point $x_j \in \X$
  has some
probability distribution $P_j$ (we write this $x_j \sim P_j$), where for $j>1$ the
distributions $P_j$ may depend on the previous points
$x_1,\ldots,x_{j-1}$ and on the results of the objective function
evaluations at these points (the function evaluations may not be
noise-free).
We say that this algorithm converges if
for any $\delta\!>\!0$, the sequence of points $x_j$ arrives at  the set
$W
(\delta)\!=\!\left\{x\in \X\colon f(x)\!-\!f_* \! \leqslant \!\delta\right\}$ with
probability one. If the objective function is evaluated without error then this obviously implies convergence (as $n \to \infty$)
 of record values  $f_{{\rm o},j}=\min_{i=1...j}  f(x_i)$
 to   $f_*$ with probability~1.

In view of continuity of $f(\cdot)$ in the neighbourhood of $\X_*$,
the event of arrival of sequence of points $x_j$  at  the set
$W
(\delta)$ with given $\delta>0$,
is equivalent to the arrival of this  sequence at  the set
$B_*(\varepsilon)= \cup_{x_* \in \X_*} B(x_*,\varepsilon)
$ for some $\varepsilon>0$ depending on~$\delta$.

Conditions on the distributions $P_j$ ($j=1,2,\ldots$)
ensuring convergence of the GRS algorithms are well understood; see, for
example, \cite{pinter,solis} and \cite[Sect.~3.2]{Z}.
Such results  are consequences of the classical in probability theory `zero-one law' or Borel-Cantelli lemmas (see e.g. \cite[Section 7.3]{grimmett2020probability}) and provide sufficient
conditions on convergence. We follow  \cite[Theorem 2.1]{zhigljavskystochastic}  to provide the most general sufficient conditions for convergence of  GRS algorithms.

{\bf Theorem 1.} {\it Consider a GRS algorithm with $x_j\sim P_j$ and let  $B \subset \X$ be a Borel subset of $\X$.
  Assume that
\be
\sum_{j=1}^\infty q_j(B)=\infty \, ,
 \label{eq:Borel1}
\ee
where $q_j(B)= \inf P_j(B)$ and the infimum is taken over all locations of previous points $x_i$ ($i=1, \ldots, j-1$) and corresponding results of evaluations of $f(\cdot)$.
Then  the sequence of points $\{x_1, x_2, \ldots\}$ falls infinitely often into the set $B$, with
probability 1. }

Note that Theorem~1 
does not make any assumptions about observations of $f(\cdot)$ and hence is valid
for the very general case where
evaluations of the objective function $f(\cdot)$ are noisy and
the noise is not  necessarily random.

Consider the following three particular cases.

(a) If in \eqref{eq:Borel1}  we use  $B=B_*(\varepsilon)$ or $B=W
(\varepsilon)$ with some  $\varepsilon>0$,
then  Theorem~1 
 gives a sufficient condition that  the corresponding GRS algorithm converges; that is,  there exists
  a subsequence
$\{x_{i_j}\}$ of the sequence $\{x_{j}\}$ which converges (with probability 1)
to the set
$\X_*$  \AZ{in the sense that 
the distance between $x_{i_j}$ and $\X_*$ tends to $0$ as $j \to \infty$.} For this subsequence $\{x_{i_j}\}$, we have $f(x_{i_j}) \to f_*$ as $j \to \infty$.

If the evaluations of $f(\cdot)$ are noise-free, then we can use the sequence of record points (that is, the points  where the records $f_{{\rm o},j}= \min_{\ell<j} f(x_\ell) $ are attained) as  $\{x_{i_j}\}$; in this case,
$f(x_{i_j})=f_{{\rm o},j}$ is the sequence of records
 converging to $f_*$ with
probability~1.
By the dominated convergence theorem (see e.g. \cite[Section 7.2]{grimmett2020probability}),
convergence of the sequence of records
$f_{{\rm o},j}$  to $f_*$ with
probability~1 implies other important types of convergence of $f_{{\rm o},j}$  to $f_*$ --- in mean and  mean square:
$
E f_{{\rm o},j} \to f_*
$
and $
E (f_{{\rm o},j}-f_*)^2 \to 0
$
as $ j \to \infty.$

(b) If \eqref{eq:Borel1} holds for  $B=B(x,\varepsilon) $   with any $x\in \X$ and any $\varepsilon>0$,
then  Theorem~1 
 gives a sufficient condition that  the sequence of points $\{x_1, x_2, \ldots\}$ is dense
 with probability 1. As this is a stronger sufficient condition than in (a), all  conclusions of (a)  are valid.

(c) If we use pure random search (PRS)  with $P=P_U$, the   uniform distribution on  $\X$ (that is, $P_j=P_U$ for all $j$ and the points $x_1, x_2, \ldots$ are independent), then the assumption that $\X$ is convex  implies $B(x,\varepsilon)\geq {\rm const}_\varepsilon>0$
for all $x\in \X$ any $\varepsilon>0$
and therefore the condition (\ref{eq:Borel1}) trivially holds for any $B=B(x,\varepsilon) $, as in (b) above.
In practice, the usual choice of the distribution $P_j$ is
\be
\label{eq:Borel4}
P_{j}=\alpha_{j}P_U+(1-\alpha_{j})Q_{j}\, ,
\ee
where $0\leqslant \alpha_{j}\leqslant 1$
and $Q_j$ is a specific  probability measure on $\X$ which may depend on
 previous evaluations of the objective function.
Sampling from the
distribution (\ref{eq:Borel4}) corresponds to taking a
uniformly distributed random point in $\X$ with
probability $\alpha_{j}$ and sampling from $Q_j$ with probability
$1-\alpha_{j}$.
In case of distributions (\ref{eq:Borel4}), the condition
$\sum_{j=1}^\infty\alpha_j=\infty$
yields the fulfilment of \eqref{eq:Borel1} for all $B=B(x,\varepsilon) $
and therefore the  GRS algorithm with such $P_j$ is theoretically converging.

\subsection{Rate of convergence}\label{sec:rate_of_convergence}


 Consider first a PRS algorithm, where  $x_j$ are i.i.d. with distribution $ P$.
 Let $\varepsilon, \delta >0$  be fixed and $B$ be the target set we want to hit by points  $x_1,x_2, \ldots$.
 For example, we set
$B=W (\delta)=\{x\in \X\colon f(x)-f_*\leqslant \delta\}$
in the case when the accuracy is expressed in terms of closeness with respect to the function value,
 $B=B(x_*,\varepsilon)$
if we are studying convergence  towards a particular global minimizer $x_*$, and $B=B_*(\varepsilon)$ if the aim is to approach a neighbourhood
of
$\X_*$.

Assume that $P$ is such that  $P(B)\!>\!0$. In particular, if $P=P_U$  is the uniform probability measure on $\X$, then, as $\X$ has  { Lipschitz} boundary, we have
  $P(B)= {\rm vol}(B) / {\rm vol}(\X) \!>\!0$.
   Note that in all interesting instances the value $p=P(B)$ is positive but small, and this will be assumed below.

Define the Bernoulli trials where the success in the trial $j$ means  $x_j \in B$.
PRS generates a sequence of
independent Bernoulli trials with the same success probability
$ {\rm Pr}\{x_j\in B\}=P(B)$.
In view of  independence of $x_1,x_2, \ldots$, we have
\bea
\textrm{Pr}\{x_1\notin B,\ldots,x_n\notin B\}=
\left(1-P(B)\right)^n
\eea
and therefore the probability
\bea
\label{eq:tends_to_zero} \textrm{Pr}\{x_j\in B\textrm{ for at
least one }j,\; 1\leqslant j\leqslant n \}=
1-\left(1-P(B)\right)^n
\eea
 tends to one as $n\to\infty$.

Let  $n_\gamma$ be the number of points which are required for PRS to reach the set
$B$ with probability at least $ 1-\gamma$, where
$\gamma\in(0,1)$; that is,
$$
n_\gamma = \min\{ n : \; 1-\left(1-P(B)\right)^n\geqslant 1-\gamma \}\, .
$$
Solving the equation $1-\left(1-P(B)\right)^n\geqslant 1-\gamma$ with respect to $n$, we obtain
\be
\label{eq:n_g}
n_\gamma= \left\lceil{\ln\gamma}/{\ln\left(1-P(B)\right)} \right\rceil\, \cong {(- \ln \gamma})/{P(B)}
\ee
as $P(B)$ is small and $\ln\left(1-P(B)\right)\cong-P(B)$ for small $P(B)$.

The numerator $- \ln \gamma$  in the expression \eqref{eq:n_g} for $n_\gamma$ depends on $\gamma$ but it is not large; for example, $-\ln\gamma \simeq 4.605$ for $\gamma=0.01$.
However, the denominator $P(B)$ (depending on  $\varepsilon$,  $d$ and the shape of $\X$)  can be very small.

Assuming that  $B=B(x_*,\varepsilon)$, where the norm is standard Euclidean, and $B$ is fully inside $\X$, we have
\be\label{eq:ball}
 {\rm vol} (B(x_*,\varepsilon) )={\rm vol} ({\cal B}(x_*,\varepsilon))=  V_d \, \varepsilon^d  \,,
\ee
where
$ \label{eq:ball0}
V_d= {\pi}^{d/2} / \left[\Gamma (d/2\!+\!1)\right]\,
$
 is the volume of the unit Euclidean ball ${\cal B}(0,1)$ and  $\Gamma(\cdot)$ is the gamma-function.
The resulting version of the expression \eqref{eq:n_g} for $n_\gamma $ in the case $B=B(x_*,\varepsilon)$ and vol$(\X)=1$ becomes
\be
\label{eq:n_gAS}
n_{\gamma}^{\rm as}
= {- \ln \gamma}/\left(\varepsilon^d V_d\right) \,.
\ee

 As $\varepsilon \to 0$, the ball $B=B(x_*,\varepsilon)$ lies fully inside $\X$ for $P_U$-almost all $x_*$.
 Indeed, asymptotically, as $n \to \infty$, the covering radius computed for uniformly distributed
random points $x_j$, tends to 0 and hence the equality \eqref{eq:ball} is valid asymptotically for almost all $x_*$. This is the reason for superscript  `as' in
\eqref{eq:n_gAS}.
   As shown below in Section~\ref{sec:non_classical},  in the non-asymptotic regime in situations (M) and (H), the volume ${\rm vol} (B(x_*,\varepsilon))$ is necessarily smaller than given by \eqref{eq:ball} and therefore the true $n_\gamma$ is (much) larger than $n_{\gamma}^{\rm as}$ in \eqref{eq:n_gAS}.

%
%


Consider now general GRS algorithms where  the probabilities $P_j$ are chosen  in the form
(\ref{eq:Borel4}), where the coefficients $\alpha_j$ satisfy the condition (\ref{eq:Borel1}).
Instead of the equality $\textrm{Pr}\{x_j\in B\}=P(B)$ for all \mbox{$j\geqslant 1$}, we now have the inequality
$\textrm{Pr}\{x_j\in B\}\geqslant \alpha_j P_U(B),$
where the equality holds in the worst-case scenario.
We define  ${n(\gamma)}$ as the smallest integer
such that the inequality
$\sum_{j=1}^{n(\gamma)} \alpha_j  \geqslant - {\ln \gamma}/{P_U(B)}\, $
is satisfied.
For the choice $\alpha_j =1/j$, which is a common recommendation, we can use the approximation $\sum_{j=1}^{n} \alpha_j \simeq \ln n $.
Therefore we obtain
$n(\gamma) \simeq \exp \{ - {\ln \gamma}/{P_U(B)} \}$.
For the case of $\X =[0,1]^d$ and $B=B(x_*,\varepsilon)$,
we obtain $n(\gamma) \simeq \exp \{ c \cdot \varepsilon^{-d} \}$, where $c =  ( - {\ln \gamma})/V_d$ with $V_d=[\sqrt{\pi}  / \left[\Gamma (d/2\!+\!1)\right]$, the volume of the unit ball.
Note also that if the distance between $x_*$ and the boundary of $\X$ is smaller than $\varepsilon$,
then the constant $c$ and hence $n(\gamma) $ are even larger.
For example, for $\gamma=0.1$, $d=10$ and $\varepsilon=0.1$,  $n(\gamma)$ is larger than $10^{ 1000000000}$.
Even for optimization problems in a small dimension $d=3$, and for $\gamma=0.1$ and $\varepsilon=0.1$,
the number $n(\gamma)$ of points required for the GRS algorithm to hit the set $B$ in the worst-case scenario is huge:  $n(\gamma) \simeq   10^{238}$.

%
%
%

\section{Points uniformly distributed on $\X$}

\label{sec:non_classical}

\label{sec:unif}

\subsection{Asymptotic case}

In this section, the point set $X_n=\{x_1, \ldots, x_n\}$ consists of the first $n$ points of a sequence  $X_\infty=\{x_1,  x_2, \ldots\}$ of independent uniformly distributed random vectors in $ \cal X$. Assume, without loss of generality, that ${\rm vol}(\X)=1$.

Consider the random variable $\rho (U,X_n)$, the distance between $U $ (the uniform random point in $\cal{X}$) and $X_n$; see \eqref{eq:CR5} for the definition of $\rho$.
The cdf (cumulative distribution function) of  $\rho(U,X_n)$ gives  the average proportions of $\X$ which are  covered by the balls centered at $X_n$ with radius $r$. That is,
\be
\label{eq:cdf_1}
F_d(r,X_n) := {\rm Pr}(\rho(U,X_n)\leq r ) ={\mathbb{E}_{X_n}{\rm vol}(B(X_n,r))} \,,
\ee
where the set $
B(X_n,r)$ is  defined in \eqref{eq:covering}. In asymptotic considerations, we need to suitably normalize the radius (which tends to zero as $n \to \infty$) in  \eqref{eq:cdf_1}. We thus consider the following sequence of cdf's:
\be
\label{eq:cdf_2}
F_n(t):={\rm Pr}( n^{1/d} V_d^{1/d} \rho(U,X_n) \leq t )= F_d\left( [nV_d]^{-1/d}\, t,X_n \right) \,,
\ee

\begin{lemma} \label{lem:1}
\be
\label{eq:Zador1}
 F_n(t) \rightarrow F(t):= 1-\exp(-t^{d})    \;\; \mbox{ as $n\rightarrow \infty$}\,.
\ee
where \AZ{the convergence is uniform in $t$} and cdf's $F_n$ are defined in \eqref{eq:cdf_2}.
\end{lemma}

The statement of Lemma~\ref{lem:1} follows from Zador's arguments in his fundamental paper  \cite{zador1982asymptotic}; see the beginning of page 142. The key observation of Zador
 is that asymptotically, as $n \to \infty$, the covering radius computed for uniformly distributed
random points $x_j$, tends to 0 and hence the equality \eqref{eq:ball} is valid asymptotically for almost all $U$; this is formula (19) in  \cite{zador1982asymptotic}.
The statement of  Lemma~\ref{lem:1} is in fact a particular case of Theorem 9.1 in \cite{graf2007foundations}, if  $Q$ is chosen as the uniform distribution on $\X$.

In what follows, we will need the  $(1-\gamma)$-quantile ($0 < \gamma<1$)
 of the cdf $F$ in the rhs of \eqref{eq:Zador1}. This $(1-\gamma)$-quantile is determined as  $t_{1-\gamma}=[-\log(\gamma)]^{1/d}$, for which we have  $F(t_{1-\gamma})=1-\gamma$ . The quantity $t_{1-\gamma}$ can be interpreted as the normalised asymptotic radius required for  covering  a subset of $\X$ of volume $(1-\gamma)$ (the weak covering introduced in Section~\ref{sec:intro}). For very small~$\varepsilon$, to cover  a subset of $\cal X$ with random centers $x_j \in X_n$ of volume which is approximately $1-\gamma$, $n=n_{\gamma}$ should satisfy
\be\label{n_gamma}
n_{\gamma} = \frac{t_{1-\gamma}^d}{\varepsilon^{d}V_d}=  \frac{-\ln(\gamma)}{\varepsilon^dV_d} \,,
\ee
which coincides with \eqref{eq:n_gAS}.
The above result can be reformulated in terms of the asymptotic radius $r$ as follows: for very large~$n$  the union of $n$ balls with random centers $x_j \in X_n$ and  radius
\be\label{normalised_radius}
r_{n,1-\gamma} = n^{-1/d}{V}_d^{-1/d}t_{1-\gamma} =  n^{-1/d}{V}_d^{-1/d}[-\log(\gamma)]^{1/d} \,
\ee
covers a subset of $\cal X$ of volume which is approximately $1-\gamma$.

In the non-asymptotic (finite $n$) regime, the distribution function $F_d(r,X_n)$ of \eqref{eq:cdf_1} can be obtained in the following way (below, for $X_n=\{x_1, \ldots, x_n\}$, the components $x_1,  x_2, \ldots x_n$ are not necessarily uniform but are i.i.d.).

\AZ{Conditionally on $U$, we have for fixed $U  \in {\cal{X}}$:} 
\be
\mathbb{P} \left\{ U \in {\cal B}_d(X_n,r)  \right\}&=& 1-\prod_{j=1}^n \mathbb{P} \left\{ U \notin {\cal B}_d({x}_j,r)  \right\} \nonumber \\
&=& 1-\prod_{j=1}^n\left(1-\mathbb{P} \left\{ U \in {\cal B}_d({x}_j,r)  \right\} \right) \nonumber\\
&=& 1-\bigg(1-\mathbb{P}_{X} \left\{ \|U - {X} \| \leq r \right\} \bigg)^n\, ,
\label{eq:prod}
\ee
where $X$ has the same distribution as $x_1$.
From \eqref{eq:prod},
the distribution function $F_d(r,X_n)$ can be obtained by averaging over {the distribution of $U$:}
\be \label{eq:prod5}
F_d(r,X_n)=  \mathbb{E}_{_U} \mathbb{P} \left\{ U \in {\cal B}_d(X_n,r)  \right\}\, .
\ee
 For  large $n$ and small $r$ we use an approximate equality $\mathbb{P}_{X} \left\{ \|U - {X} \| \leq r \right\} \simeq r^{d}V_d $ in \eqref{eq:prod}. By doing so, averaging with respect to $U$ is redundant and we arrive at the results of Section~\ref{sec:rate_of_convergence}. If $n$ is not so large, the quantity $\mathbb{P}_{X} \left\{ \|U - {X} \| \leq r \right\}$ has to be approximated by other means. This will be discussed in Section~\ref{Petrov_section}.


\subsection{Bounds for  $F_d(r,X_n)$}\label{sec:jensen}
Evaluating the expectation in \eqref{eq:prod5} is difficult but simple bounds can be obtained by applying Jensen's inequality. Here we will focus attention to the case of $\X=[0,1]^d$ and $X_n=\{x_1, \ldots, x_n\}$, where  $x_1,  x_2, \ldots$ is a sequence of uniformly distributed random vectors on $\X$. From \eqref{eq:prod5}, we have
\bea
 \mathbb{E}_{_U} \mathbb{P} \left\{ U \in {\cal B}_d(X_n,r) \right\} =  1- \mathbb{E}_{_U}  \left [  \left(1-\mathbb{P}_{X} \left\{ \|U - {X} \| \leq r \right\} \right)^n  \right ]\, .
\eea
An immediate use of Jensen's inequality yields the bound:
\be \label{Jensen_naive}
 \mathbb{E}_{_U} \mathbb{P} \left\{ U \in {\cal B}_d(X_n,r) \right\} \leq  1-   \left(1-\mathbb{P}_{X} \left\{ \|  \bm{1/2} - {X} \| \leq r \right\} \right)^n  \, .
\ee
Here and below ${ \bf a}=(a,a,\ldots, a) \in \mathbb{R}^d$ for any $a$.
However, noticing the fact $\mathbb{P}_{X} \left\{ \|U - {X} \| \leq r \right\} = \mathbb{P}_{Z} \left\{ \|Z - {X} \| \leq r \right\}$ where $Z$ in a uniform random vector on $[1/2,1]^d$, we can apply Jensen's inequality to obtain:
\be \label{Jensen_main}
 \mathbb{E}_{_U} \mathbb{P} \left\{ U \in {\cal B}_d(X_n,r) \right\} \leq  1-   \left(1-\mathbb{P}_{X} \left\{ \|  \bm{3/4}-{X} \| \leq r \right\} \right)^n  \, .
\ee
The forms of  the bounds in \eqref{Jensen_naive} and \eqref{Jensen_main} suggest an approximation of the following form may be useful:
\be\label{Jensen_approx}
 \mathbb{E}_{_U} \mathbb{P} \left\{ U \in {\cal B}_d(X_n,r) \right\} \simeq  1-   \left(1-\mathbb{P}_{U,X} \left\{ \|  U-{X} \| \leq r \right\} \right)^n  \, .
\ee
Here, instead of  fixing $U$ to $\bm{1/2}$ or ${\bm{3/4}}$, it is a uniform random vector on $[0,1]^d$. The probability
${\mathbb{P}_{U,X} \left\{ \|  U-{X} \| \leq r \right\}}$
has the interpretation of being the average intersection a ball of radius $r$ with a random center at $U$ has with the cube $[0,1]^d$. For different $d$ and $r$, the distribution of ${\mathbb{P}_{X} \left\{ \|  U-{X} \| \leq r \right\}}$ normalised by the volume of the ball $r^dV_d$ is shown in Figures~\ref{key_figure9}-\ref{key_figure10}.

\subsection{Numerical studies}
In this section, we will demonstrate one of the key messages of the paper saying that in high dimensions, the asymptotic results are not attainable for reasonable values of $n$ and consequently produce poor approximations for $n$ not astronomically large.

\begin{figure}[!h]
\centering
\begin{minipage}{.5\textwidth}
  \includegraphics[width=1\linewidth]{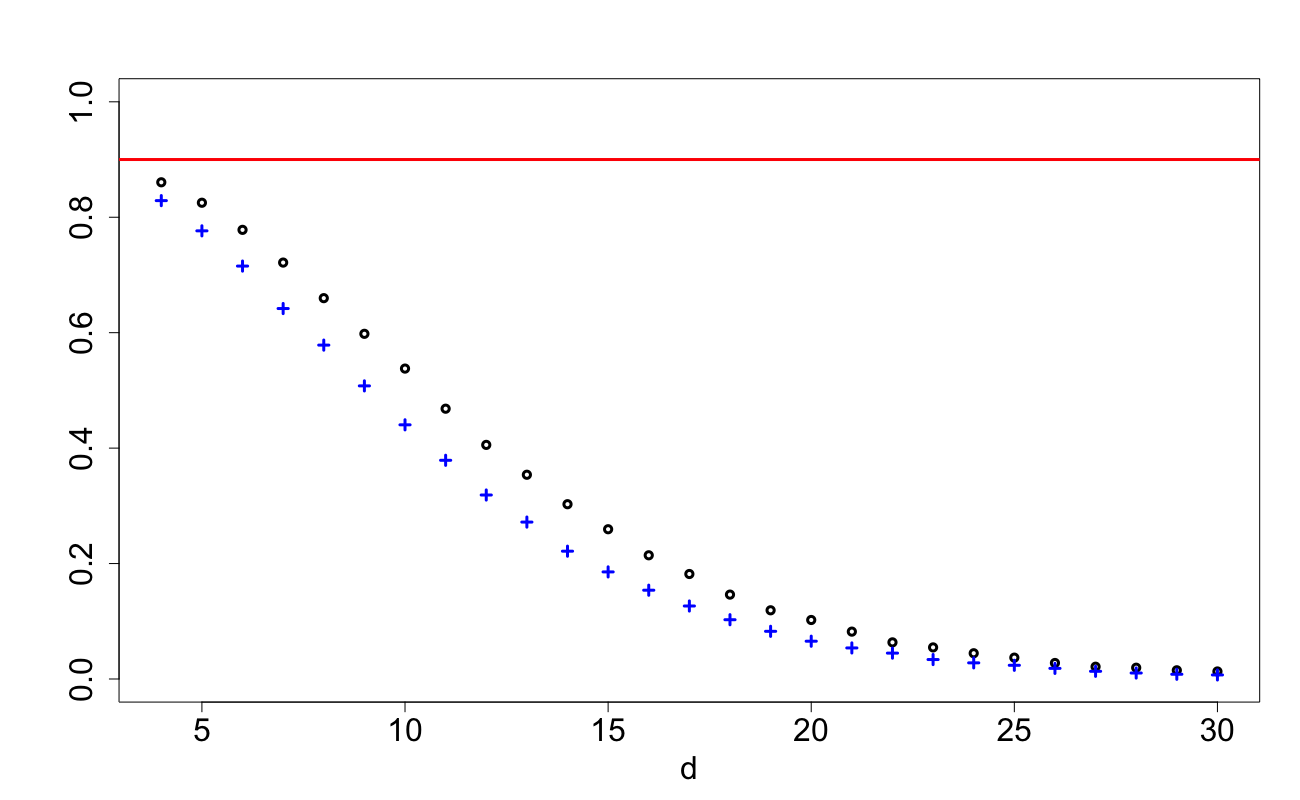}
  \caption{Covering proportions using the \\asymptotic radius; $n=1000,\;10000$. }
  \label{key_figure1}
\end{minipage}%
\begin{minipage}{.5\textwidth}
  \centering
  \includegraphics[width=1\linewidth]{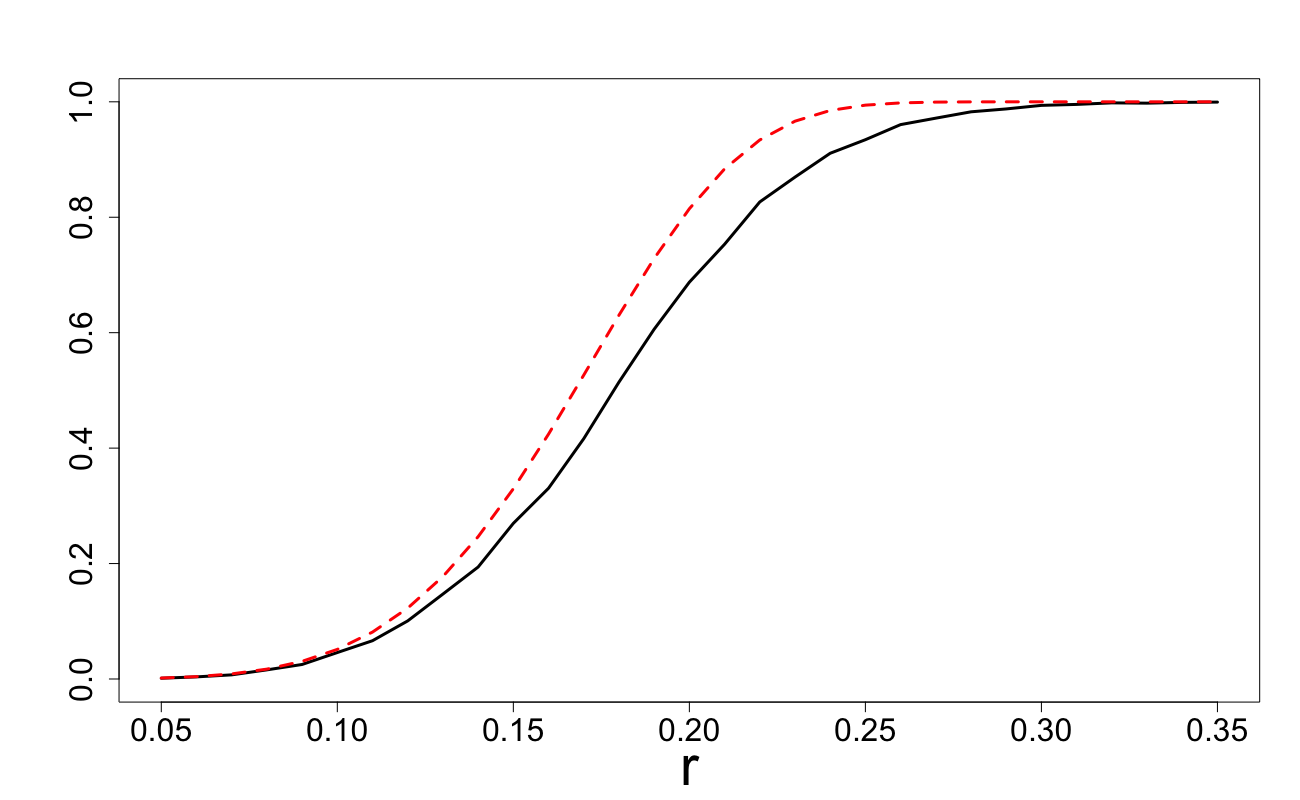}
  \caption{$ F_d(r,X_n) $ and $F(n^{1/d}V_d^{1/d}r)$ as functions of $r$; $ d=5$  and $n=1000$ }
    \label{key_figure2}
\end{minipage}
\end{figure}

In Figure~\ref{key_figure1}, we plot $F_d(r,X_n)$ as a function of $d$ for $n=1000$ (using blue plusses) and $n=10000$ (using black circles). For each value of $d$, the radius $r$ is chosen based on the asymptotic result given in \eqref{normalised_radius} with $1-\gamma=0.9$; this is shown by the solid red line at $0.9$. We see that very quickly and for $n$ that would be deemed large, $F_d(r,X_n)$ is significantly smaller than  0.9 and quickly tends to zero in $d$.

\begin{figure}[!h]
\centering
\begin{minipage}{.5\textwidth}
  \includegraphics[width=1\linewidth]{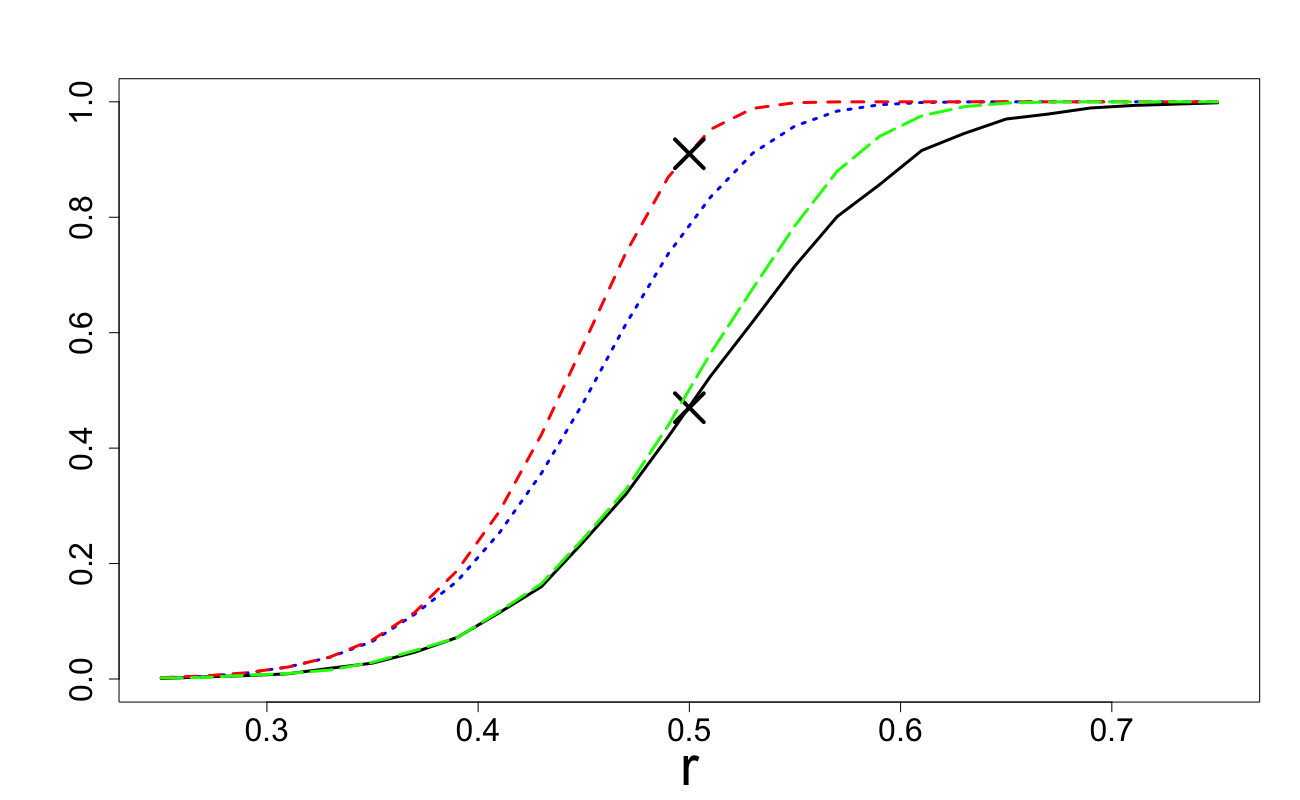}
  \caption{$d=10,\;n=1000.$ }
  \label{key_figure3}
\end{minipage}%
\begin{minipage}{.5\textwidth}
  \centering
  \includegraphics[width=1\linewidth]{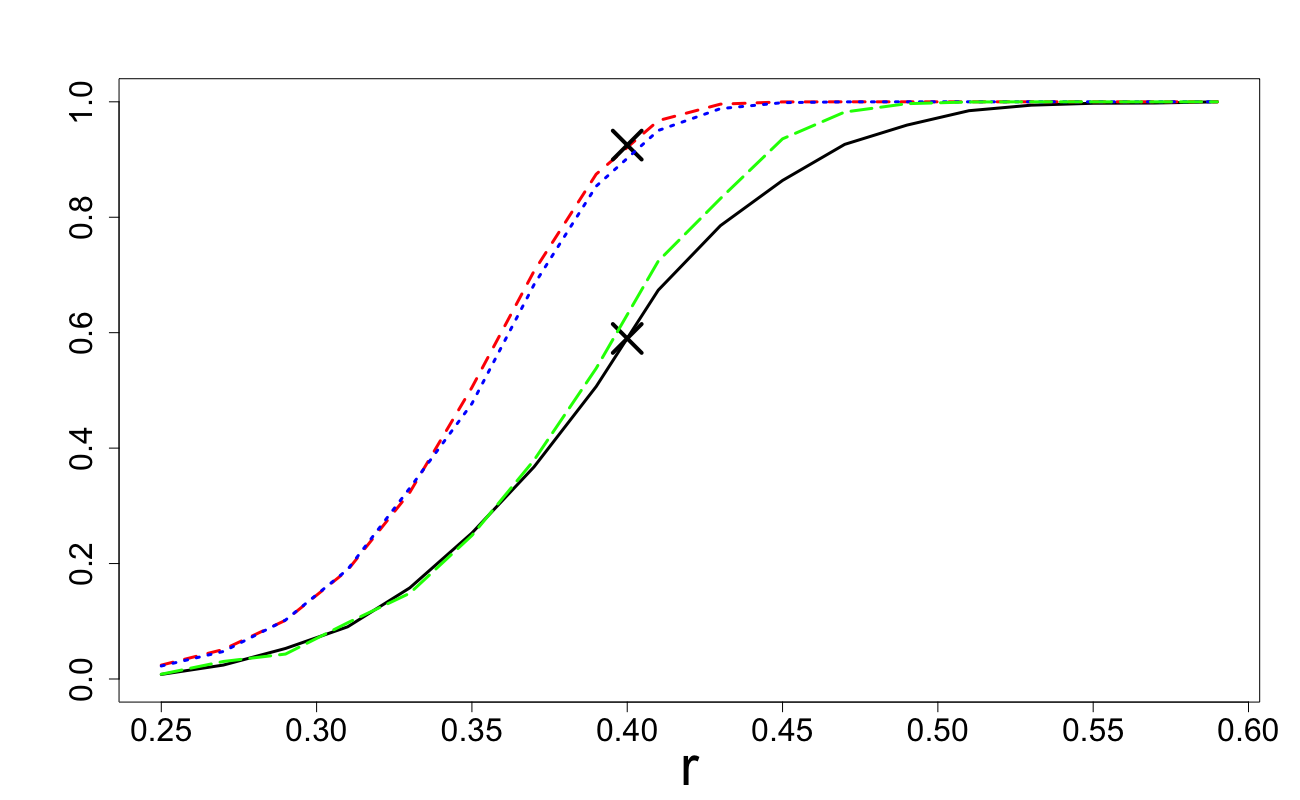}
  \caption{$d=10,\;n=10000.$ }
    \label{key_figure4}
\end{minipage}
\end{figure}

\begin{figure}[!h]
\centering
\begin{minipage}{.5\textwidth}
  \includegraphics[width=1\linewidth]{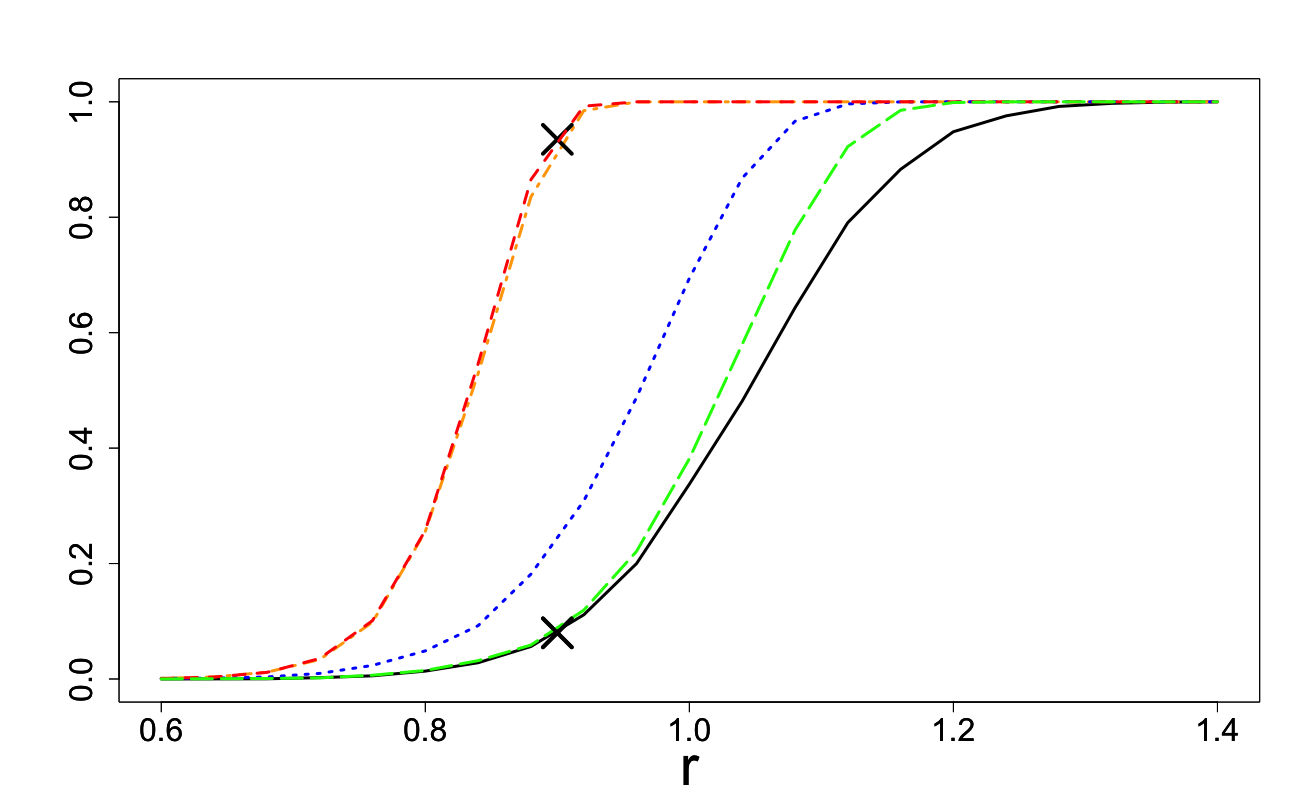}
  \caption{$d=20,\;n=1000. $ }
  \label{key_figure5}
\end{minipage}%
\begin{minipage}{.5\textwidth}
  \centering
  \includegraphics[width=1\linewidth]{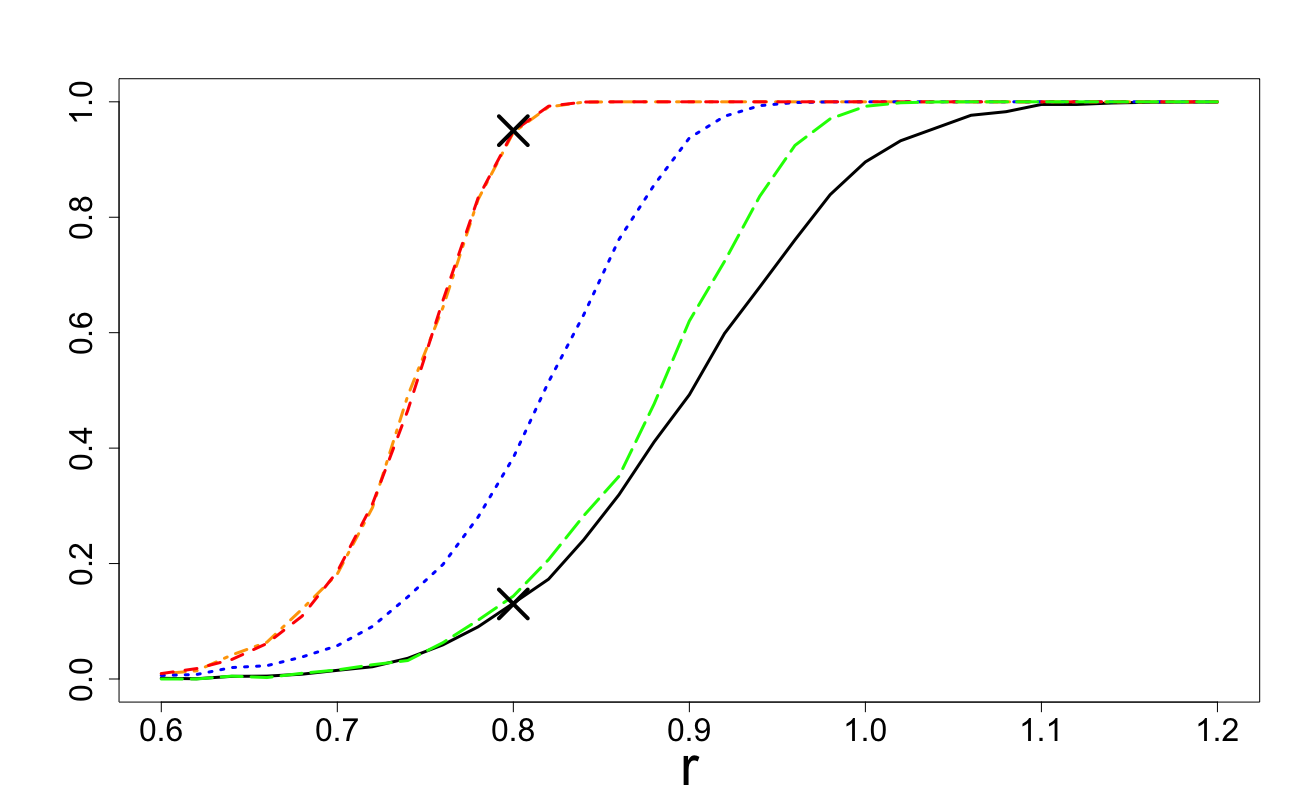}
  \caption{$d=20$, $n=10000. $ }
    \label{key_figure6}
\end{minipage}
\end{figure}

\begin{figure}[!h]
\centering
\begin{minipage}{.5\textwidth}
  \includegraphics[width=1\linewidth]{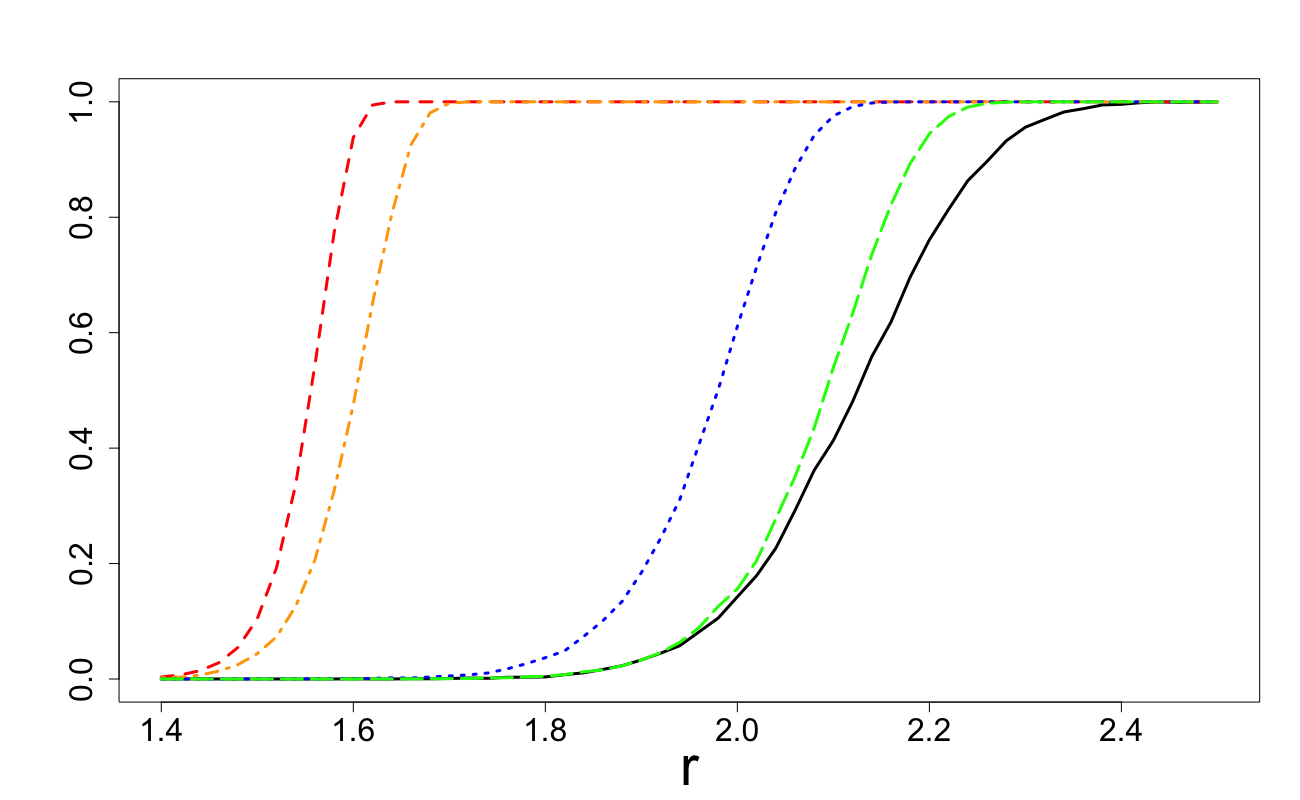}
  \caption{$n=1000, d=50$ }
  \label{key_figure7}
\end{minipage}%
\begin{minipage}{.5\textwidth}
  \centering
  \includegraphics[width=1\linewidth]{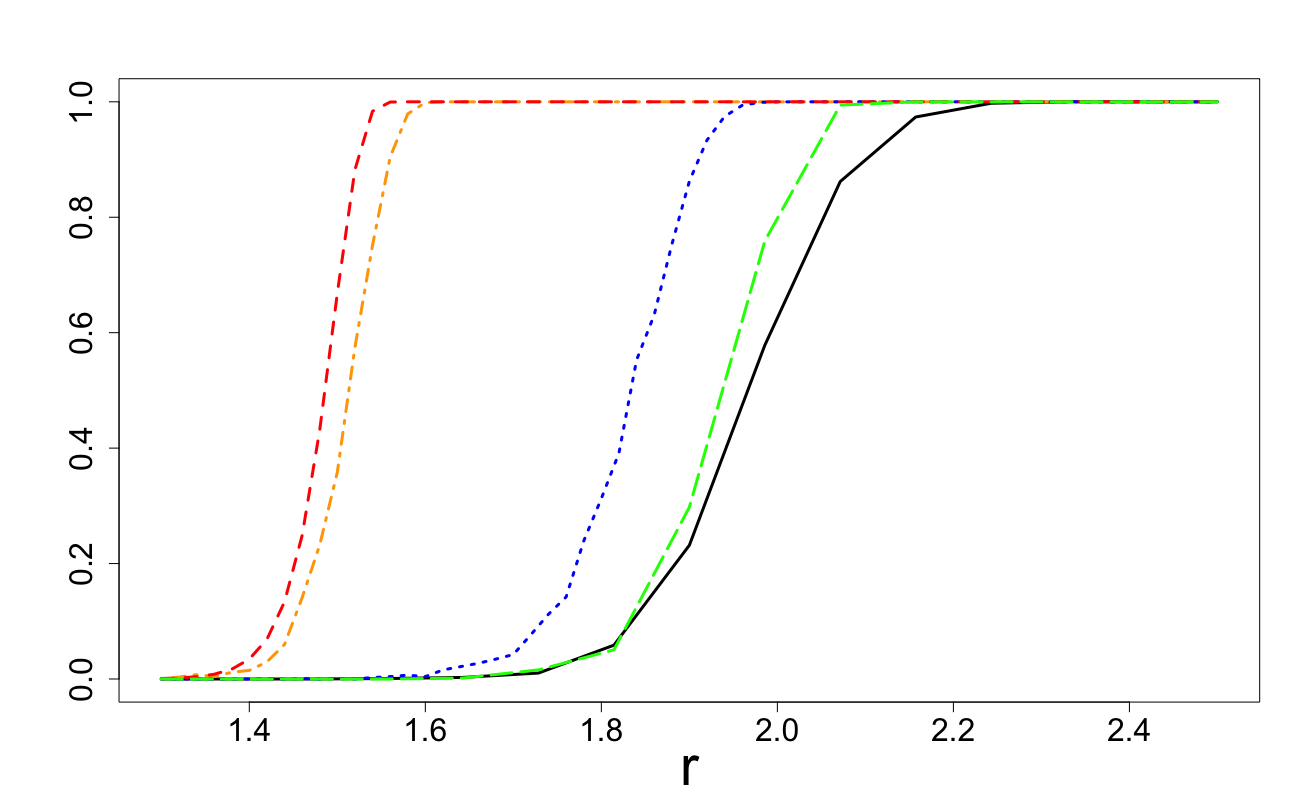}
  \caption{$n=10,000, d=50$ }
    \label{key_figure8}
\end{minipage}
\end{figure}

The big difference between the asymptotic and finite  regime is further illustrated in Figures~\ref{key_figure2}--\ref{key_figure8}. In these figures, using a solid black line we depict $F_d(r,X_n)$ as a function of $r$ for different values of $d$ and $n$ that are provided in the caption of each figure. In these figures, the dashed red line is the approximation obtained from the asymptotic result \eqref{eq:Zador1}, that is, the approximation $F_d(r,X_n) \approx F(n^{1/d}V_d^{1/d}r)$. In Figures~\ref{key_figure3}-\ref{key_figure8}, we also include two Jensen's bounds given in \eqref{Jensen_naive} (dot dashed orange) and \eqref{Jensen_main} (dotted blue), as well as the approximation given in \eqref{Jensen_approx} (longer dashed green). From these figures, we can make the following observations.
\begin{enumerate}
\item Unless $d$ is small,  the asymptotic results produce poor approximations even if $n$ is  reasonably large.
\item The approximation in \eqref{Jensen_approx} is rather accurate but worsens for smaller $\gamma$.
\item For $r\leq 1/2$, the asymptotic bounds and \eqref{Jensen_naive} coincide; this follows from the equality\\  ${\mathbb{P}_{X} \left\{ \|  \bm{1/2} - {X} \| \leq r \right\}=r^dV_d}$ for $r\leq 1/2$.
\item The refined Jensen's bound given in \eqref{Jensen_main} is  superior to \eqref{Jensen_naive} and especially to the asymptotic bound. This becomes particularly evident in higher dimensions; see Figures \ref{key_figure7} and \ref{key_figure8}.
\end{enumerate}

In Figures~\ref{key_figure3}-\ref{key_figure6}, the crosses on the dashed red line and solid black line mark points of interest. In Figure~\ref{key_figure3}, for $r=0.5$ we obtain $F(n^{1/d}V_d^{1/d}r)=0.91$ but the true value of $F_d(0.5,X_n)$ is closer to 0.41. As $n$ increases from 1000 to 10000 as is shown in Figure~\ref{key_figure4},  for $r=0.4$ we have $F(n^{1/d}V_d^{1/d}r)=0.925$ and  $F_d(0.4,X_n)$ is closer to 0.6. (Recall that in view of \eqref{eq:Zador1}, we should have $F(n^{1/d}V_d^{1/d}r) \simeq F_d(r,X_n))$ for all $r$ and $n$ large enough). The respective triples
$(r; F(n^{1/d}V_d^{1/d}r), F_d(r,X_n))$ for Figures~\ref{key_figure4}-\ref{key_figure6} are $(0.9; 0.935,0.08)$ and $(0.8; 0.95,0.13)$. For the case of $d=50$ and shown in Figures~\ref{key_figure7}-\ref{key_figure8}, the asymptotic properties are so far from being achieved with $n=1000$ and $n=10000$ that such a comparison does not even make sense.

\begin{figure}[!h]
\centering
\begin{minipage}{.5\textwidth}
  \includegraphics[width=1\linewidth]{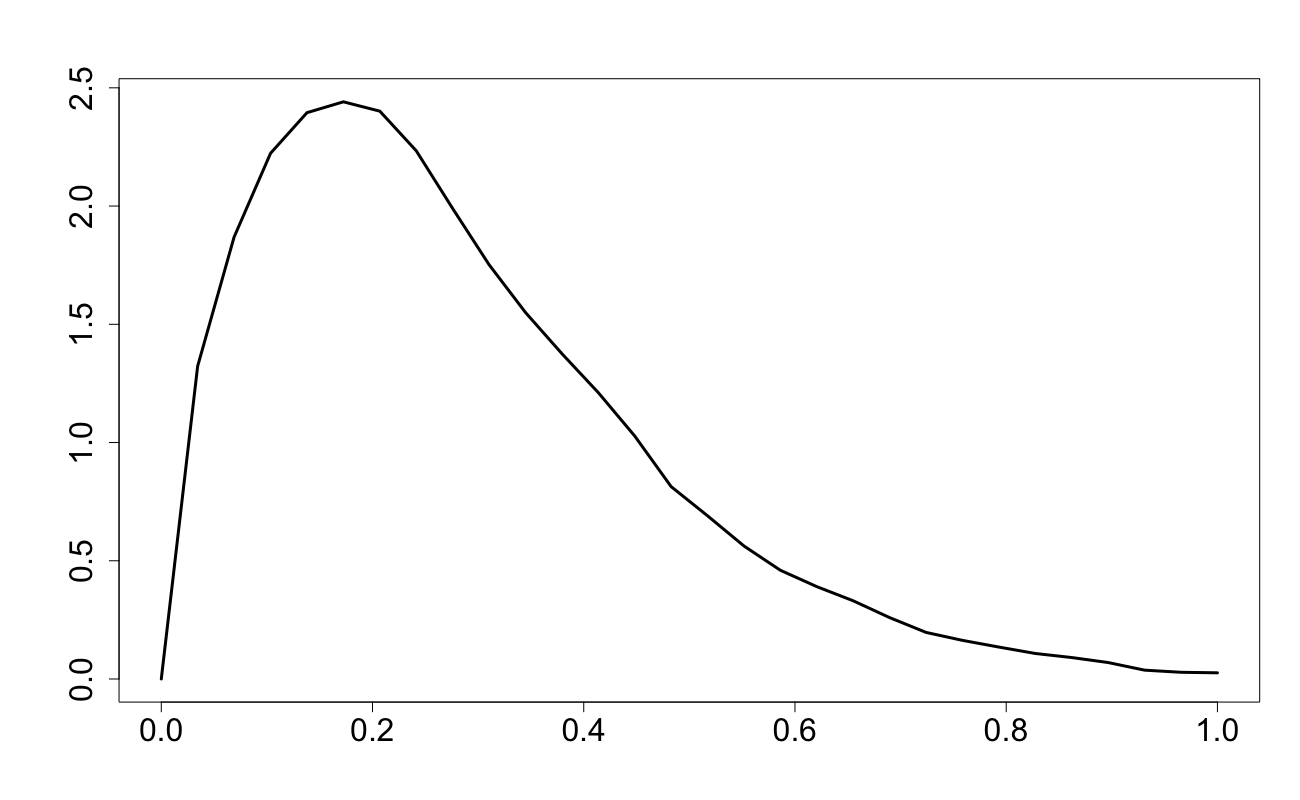}
  \caption{Density of  r.v. $\kappa_U$; $d=10,r=0.5$ }
  \label{key_figure9}
\end{minipage}%
\begin{minipage}{.5\textwidth}
  \centering
  \includegraphics[width=1\linewidth]{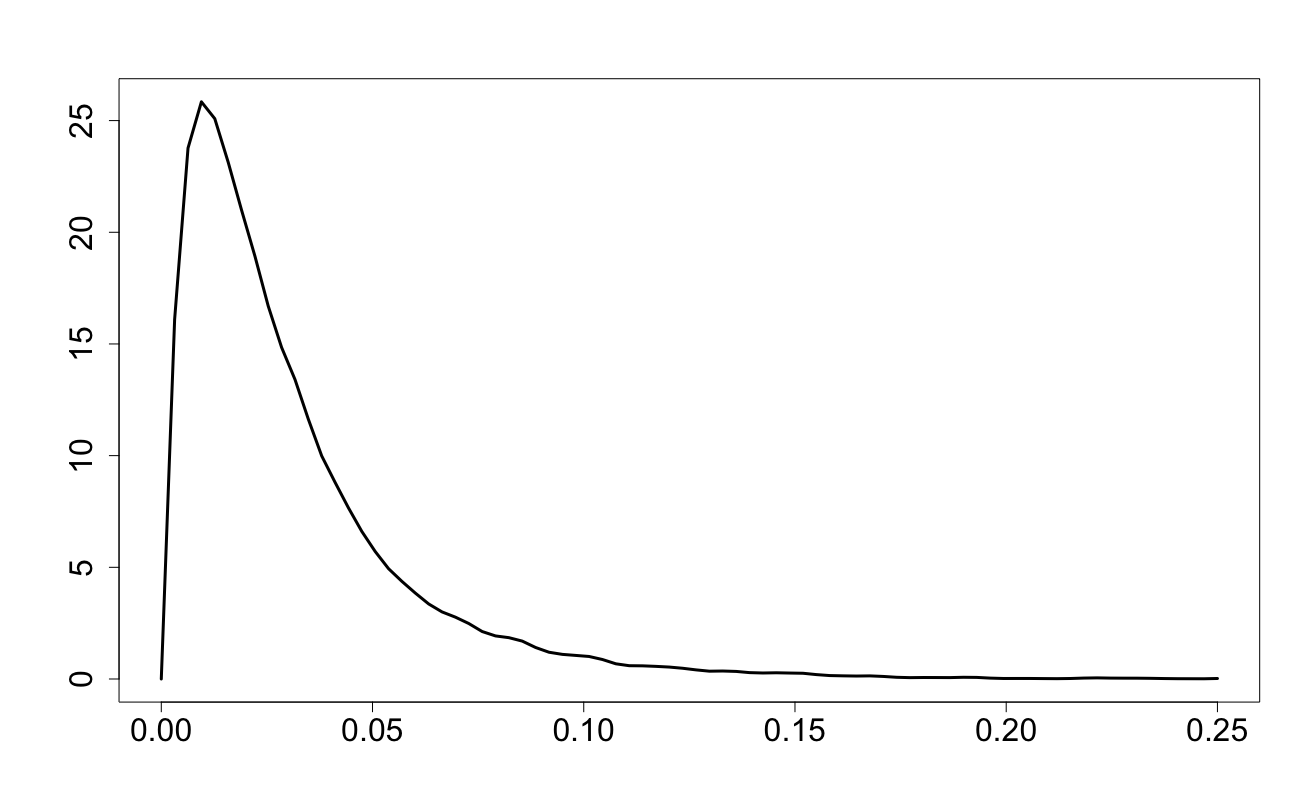}
  \caption{Density of  r.v. $\kappa_U$; $d=20,r=0.9$ }
    \label{key_figure10}
\end{minipage}
\end{figure}

In Figures~\ref{key_figure9} and \ref{key_figure10} we use $d=10$, $d=20$ and the values of $r$ corresponding to the crosses in  Figures~\ref{key_figure3} and \ref{key_figure5}. In these figures, we depict the distribution of intersection a random point has with the cube normalised by the volume of the ball $r^dV_d$; that is, we plot the density of the r.v. $\kappa_U=\mathbb{P}_{X} \left\{ \|  U-{X} \| \leq r \right\}/(r^dV_d)$, where both  $U$ and $X$ have uniform distribution on $[0,1]^d$. The importance of these two figures is another illustration of  inadequacy  of the key assumption behind \eqref{eq:n_gAS}, which can be formulated as  the assumption that the distribution of density of the r.v. $\kappa_U $ is very close to the delta-measure concentrated at one. This assumption is indeed reasonably adequate if $r$ can be chosen small enough. However, as Figure~\ref{key_figure9} and especially Figure~\ref{key_figure10} illustrate, even for relatively large values of $n$ the required values  of $r$ are not small enough for this to hold even approximately.
Note that in the derivation of the asymptotic values of $n_\gamma=n_{\gamma}^{\rm as}$ in
\eqref{eq:n_gAS} we use the value 1 rather than the random variables $\kappa_U$ with the densities shown in
Figures~\ref{key_figure9},\ref{key_figure10}.

%
%
%
%
%
%

\section{Modification of sampling schemes and non-uniform distribution of the target}

\label{sec:delta1}

\AZ{In Section~\ref{sec:unif}, we have used the principal sampling scheme where points  $x_j$ in  $X_n=\{x_1, \ldots, x_n\}$ are i.i.d. uniform on $\X=[0,1]^d$. 
In Section~\ref{sec:delta1_1} we study a modification  of this scheme where  $x_j \in X_n$
are i.i.d. uniform random points in a smaller $\delta$-cube $C_\delta=[1/2-\delta/2, 1/2+\delta/2]^d$
with $0<\delta<1$.  
In Section~\ref{sec:delta1_2} we investigate the effect of replacing random points by points from a low-discrepancy sequence. The choice of a specific low-discrepancy sequence has very little impact on the results and we present the results for Sobol sequence only.
In Section~\ref{sec:delta1_3} we  will investigate the effect of replacement of the uniform distribution of the target $x_* \in [0,1]^d$ by a bowl-shaped distribution such as the product of arcsine distributions on $[0,1]$.
}

\subsection{\AZ{Points $x_j$ are i.i.d. uniformly distributed on  $C_\delta$}  
}
\label{sec:delta1_1}

\AZ{In this section we demonstrate the $\delta$-effect, which manifests that in high dimensions sampling in a
 cube $C_\delta$ with suitable $0<\delta<1$ leads to a much more efficient covering scheme than sampling within the whole cube   $[0,1]^d$. Note that the $\delta$-effect is not  obvious  being  completely unknown in the literature on stochastic global optimization  and perhaps in literature on global optimization in general. All existing literature  recommends space-filling in the whole set $X$ and not in its subset. Moreover,  there are recommendations in literature (see, for example, \cite{janson1986random,tsvetkov2022pure})
of choosing more points closer to the boundary of the cube rather than purely uniformly in order to improve space-filling properties of  random points. }

\begin{figure}[!h]
\centering
\begin{minipage}{.5\textwidth}
  \includegraphics[width=1\linewidth]{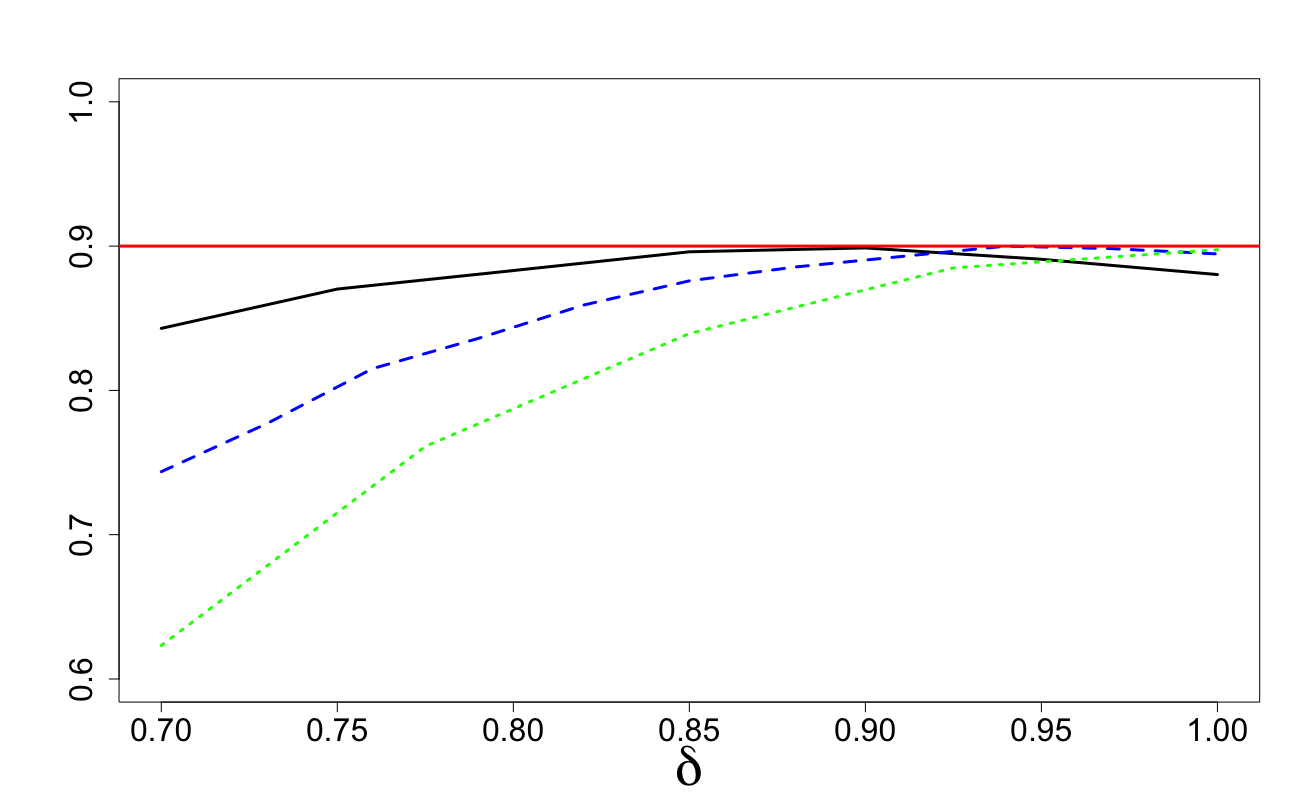}
  \caption{$d=10;\; n=1000, 10000, 100000$. }
  \label{key_figure13}
\end{minipage}%
\begin{minipage}{.5\textwidth}
  \centering
  \includegraphics[width=1\linewidth]{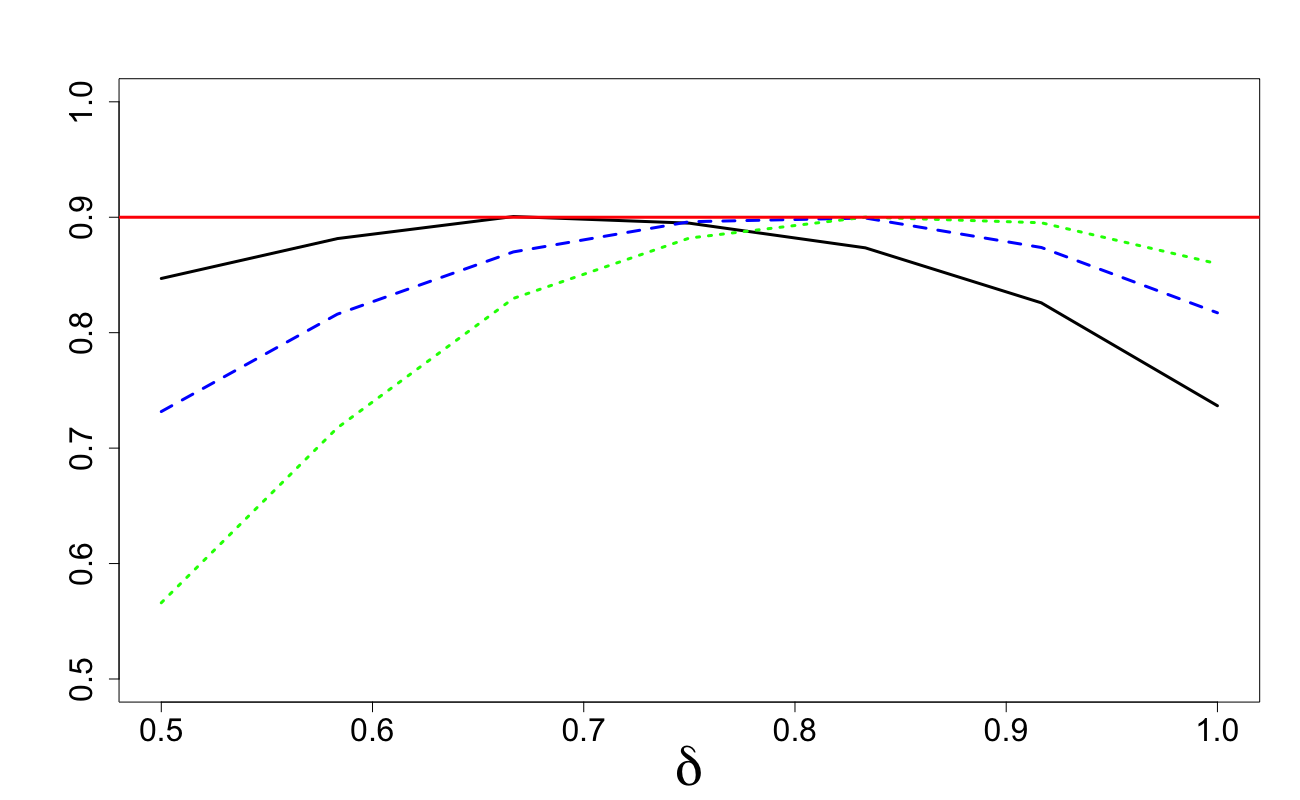}
  \caption{$d=20;\; n=1000, 10000, 100000$. }
    \label{key_figure14}
\end{minipage}
\end{figure}

\begin{figure}[!h]
\centering
\begin{minipage}{.5\textwidth}
  \includegraphics[width=1\linewidth]{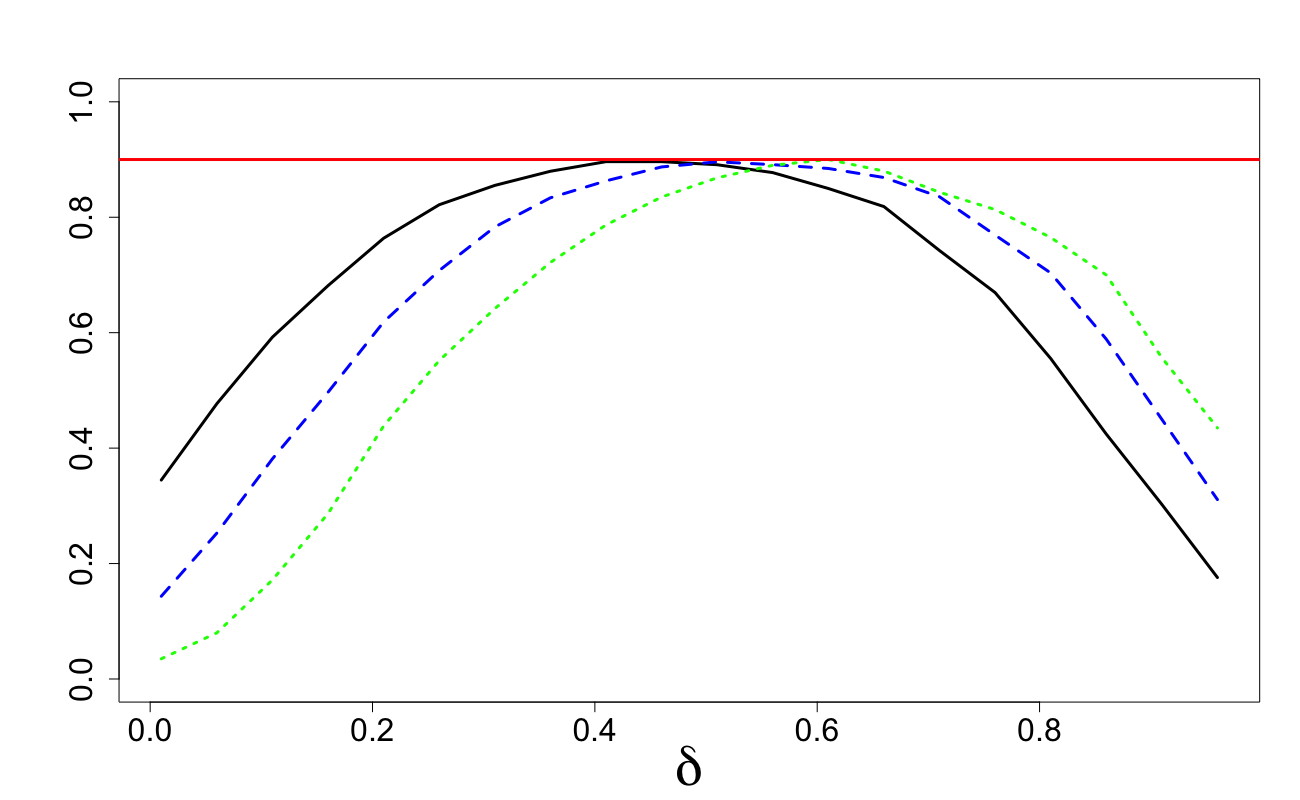}
  \caption{$d=50;\; n=1000, 10000, 100000$. }
  \label{key_figure15}
\end{minipage}%
\begin{minipage}{.5\textwidth}
  \centering
  \includegraphics[width=1\linewidth]{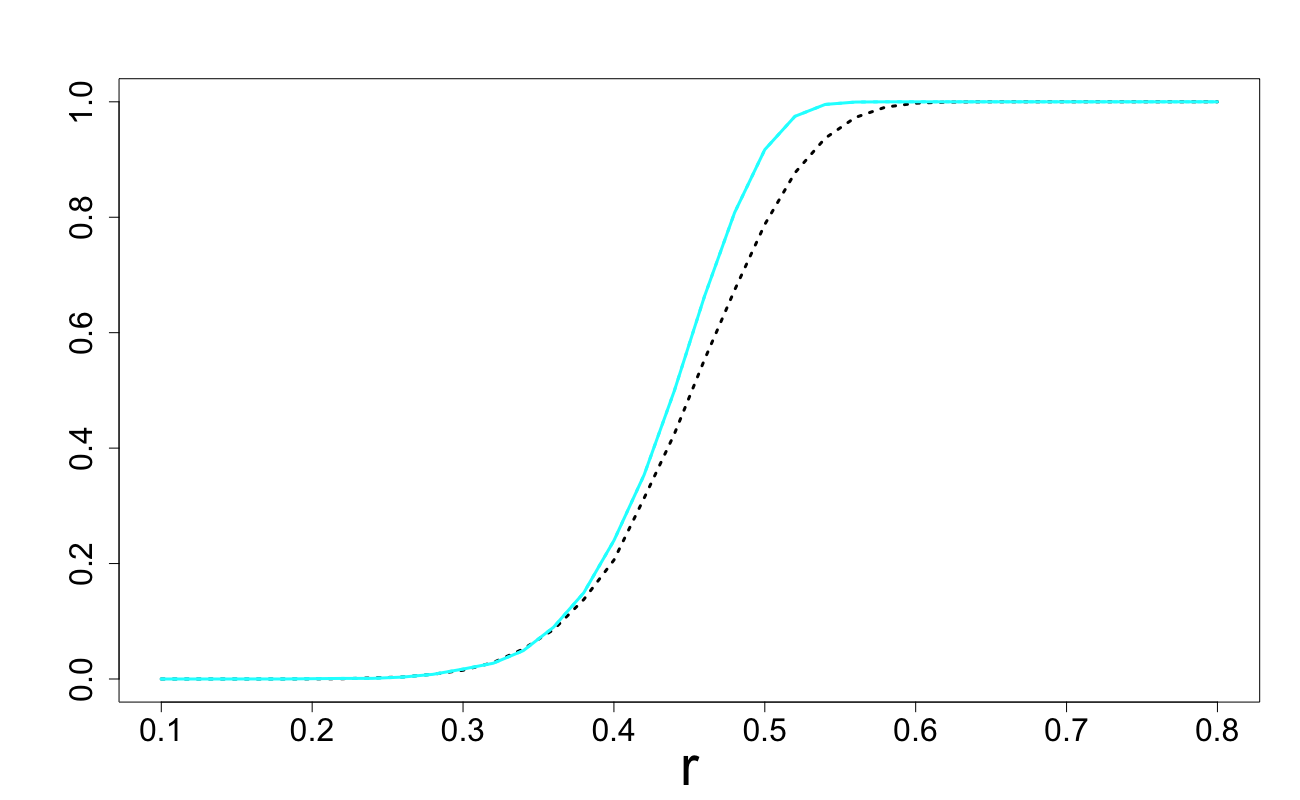}
  \caption{$d=10, n=1000$: Jensen's bound with $\delta=1$ and $\delta=0.5$. }
    \label{key_figure16}
\end{minipage}
\end{figure}

In Figures~\ref{key_figure13}--\ref{key_figure15}, for different values of $d$ and $n$ we plot $F_d(r,X_n)$ as a function of $\delta$. For each $d$ and $n$, the value of $r$ has been chosen such that $\max_{0\leq\delta\leq1}F_d(r,X_n)=0.9$; these values of $r$ (along with optimal values of  $\delta$, in brackets) can be obtained from Table~\ref{covering_table}. In these figures, the values of $F_d(r,X_n)$ for $n=1000, 10000, 100000$ are shown with a solid black line, dashed blue line and dotted green line respectively. \AZ{These figures demonstrate the  `$\delta$-effect' formulated as the second main message in Introduction. 
These figures also clearly demonstrate that despite  sampling uniformly in the cube $[0,1]^d$ is asymptotically optimal, for large $d$ it is always a poor strategy, which  can be substantially improved.}

The discussion of Jensen's bounds given in Section~\ref{sec:jensen} still apply to the case of $X_n$ sampled uniformly within $\delta $-cube $C_\delta$. The only adjustment that needs to be made to the results  of Section~\ref{sec:jensen}  is to let $X$ be a uniform random vector in $C_\delta$ and not $[0,1]^d$. In Figure~\ref{key_figure16}, we depict the Jensen's lower bound given in \eqref{Jensen_main} for $X$ uniform in $[0,1]^d$ and for $X$ uniform in the $\delta$ cube $[{1/4,3/4}]^d$ (so that $\delta=0.5$). We see that  the lower bound for $X_n$ sampled within the $\delta $-cube is larger than $X_n$ sampled from the whole cube. This further supports the conclusion that for $n$ not astronomically large, the `$\delta$-effect' should always be considered.

In Table~\ref{covering_table}, for $X_n$ chosen uniformly in the cube $[0,1]^d$ and $X_n$ chosen uniformly in the $\delta$-cube, we tabulate the values of $r_{n,1-\gamma}$ with $\gamma=0.1$ for different $d$ and $n$. In the columns labeled $\delta$-cube, the values in the brackets correspond to the approximately optimal values of $\delta$. We can see that for small $d$, the $\delta$-effect is very small (since  $n$ is relatively large in these dimensions). For larger dimensions, the $\delta$-effect is very prominent.

\begin{table}[h]
\centering
\begin{tabular}{ |p{1.2cm}||p{1cm}|p{1.5cm}|p{1cm}|p{1.5cm}|p{1cm}| p{1.5cm}|p{1cm}|p{1.7cm}| }
     \hline
 & \multicolumn{2}{c|}{$n=100$} &\multicolumn{2}{c|}{$n=1,000$} &\multicolumn{2}{c|}{$n=10,000$} &\multicolumn{2}{c|}{$n=100,000$} \\
 \hline
  & $[0,1]^d$ & $\delta$-cube& $[0,1]^d$ & $\delta$-cube& $[0,1]^d$& $\delta$-cube&  $[0,1]^d$& $\delta$-cube\\
 \hline
$d=5$ & 0.41   & 0.40 (0.9)  &   0.24   & 0.24 (1.0)  & 0.24  &  0.29 (1.0) & 0.09   &  0.089 (1.0) \\
$d=10$&  0.81 &  0.78 (0.7)&  0.61  &   0.60 (0.9) & 0.46 & 0.46 (1.0) &  0.36     &  0.360 (1.0)  \\
$d=15$& 1.13 & 1.04 (0.6)  &  0.91 & 0.88 (0.8)   & 0.76  & 0.74 (0.9)  &   0.62   &  0.619 (0.9)  \\
$d=20$& 1.38   & 1.25 (0.5)  &  1.17& 1.11 (0.7)  & 1.01 & 0.97 (0.8) &  0.87     & 0.855 (0.9)    \\
$d=25$&  1.60 & 1.42 (0.5) &  1.39 & 1.30 (0.6)   &1.23  & 1.18 (0.8) &    1.09   &  1.060 (0.8)  \\ 
$d=50$& 2.46  & 2.07 (0.4) &  2.26 &  1.98 (0.5)  & 2.10 & 1.90 (0.5) & 1.96    & 1.790  (0.6)  \\
 \hline
\end{tabular}
\caption{Values for $r_{n,1-\gamma}$ with $\gamma=0.1$.}
\label{covering_table}
\end{table}

In Tables~\ref{n_gamma_1}--\ref{n_gamma_3}, we consider an equivalent reformulation of the results of Table~\ref{covering_table}. In these tables, for a given $r$ we specify the value of $n_\gamma$, with $\gamma=0.1$, for $X_n$ chosen uniformly in the cube $[0,1]^d$ and $X_n$ chosen uniformly in the $\delta$-cube. We also include the approximation based on the the asymptotic arguments leading to \eqref{eq:n_gAS}. We see that in high dimensions, the requirement of $r$ being small enough for \eqref{eq:n_gAS} to provide  sensible approximations requires $n$ to be extremely large. Such large values of $n$ are impractical.

\begin{table}[h]
\centering
\begin{tabular}{|c||c|c|c|c|c|c|c|c|}
  \hline
  $r$ & 0.9 & 0.95 & 1 & 1.05 & 1.1 & 1.15 \\
  \hline
  $n_\gamma$ with $\delta=1$ & 54,000 & 25,000 & 10,800 &4,600  &2,700  & 1,300 \\
  $n_\gamma$ with $\delta=\delta^*$ & 40,000 (0.8) & 15,000 (0.8) & 6,300 (0.8) & 2,700 (0.7) & 1,100 (0.7)  & 500 (0.7) \\
   $n_\gamma$ from \eqref{eq:n_gAS}  &734 &  249 & 89 & 34 & 13 &5  \\
  \hline
\end{tabular}
\caption{Values of $n_{\gamma}$: $d=20, \gamma=0.1.$}
\label{n_gamma_1}
\end{table}


\begin{table}[h]
\centering
\begin{tabular}{|c||c|c|c|c|c|c|c|c|c|}
  \hline
  $r$ & 2& 2.05 & 2.1 & 2.15 & 2.2 & 2.25& 2.3 \\
  \hline
  $n_\gamma$ with $\delta=1$ & 50,000 & 21,000 & 10,000 &5,000  &2,200  & 1,200& 600\\
  $n_\gamma$ with $\delta=\delta^*$ &700 (0.4) & 200 (0.4)  & 50 (0.3) & 12 (0.2) & 2 (0.1)  & NA & NA \\
  $n_\gamma$ from \eqref{eq:n_gAS} &0 &  0 & 0 & 0 & 0 &0&0  \\
  \hline
\end{tabular}
\caption{Values of $n_{\gamma}$: $d=50, \gamma=0.1.$}
\label{n_gamma_3}
\end{table}

\subsection{\AZ{Points $x_j$ are taken from a low-discrepancy sequence
}}
\label{sec:delta1_2}

\begin{figure}[!h]
\centering
\begin{minipage}{.5\textwidth}
  \includegraphics[width=1\linewidth]{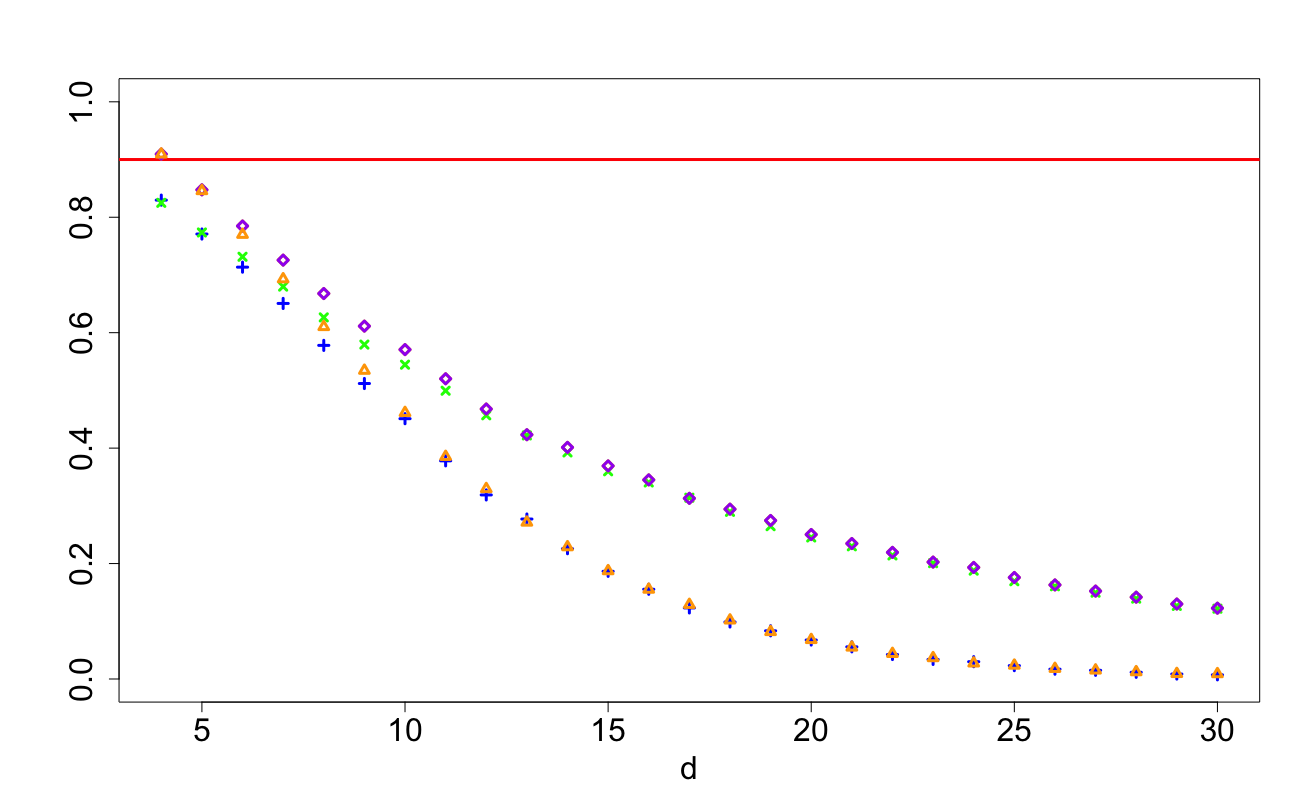}
  \caption{overing using the asymptotic radius\\ with Sobol and $\delta$-cube points: $n=2^{10}$. }
  \label{key_figure17}
\end{minipage}%
\begin{minipage}{.5\textwidth}
  \centering
  \includegraphics[width=1\linewidth]{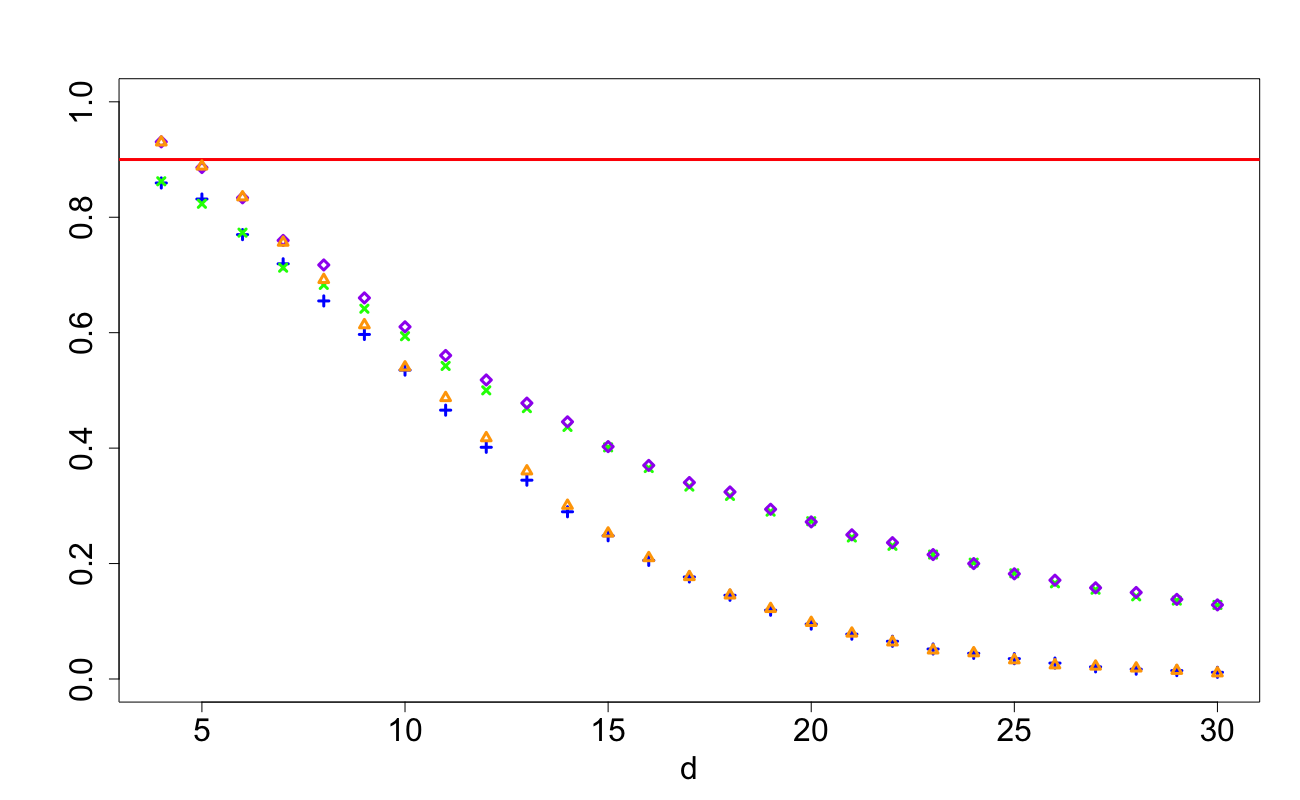}
  \caption{Covering using the asymptotic radius\\ with Sobol and $\delta$-cube points: $n=2^{13}$. }
    \label{key_figure18}
\end{minipage}
\end{figure}

\begin{figure}[!h]
\centering
\begin{minipage}{.5\textwidth}
  \includegraphics[width=1\linewidth]{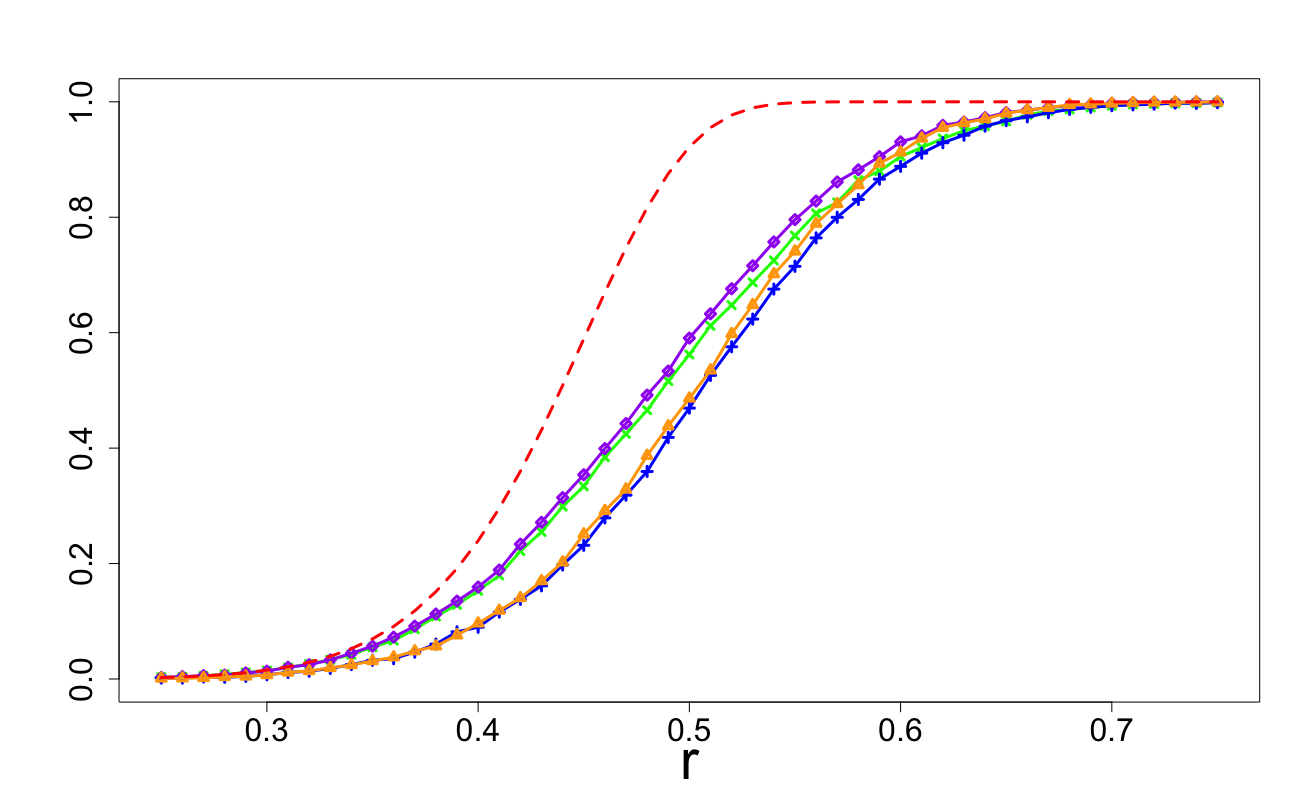}
  \caption{Sobol and $\delta$-cube points versus \\ the asymptotic covering:  $d=10, n=1024$. }
  \label{key_figure19}
\end{minipage}%
\begin{minipage}{.5\textwidth}
  \centering
  \includegraphics[width=1\linewidth]{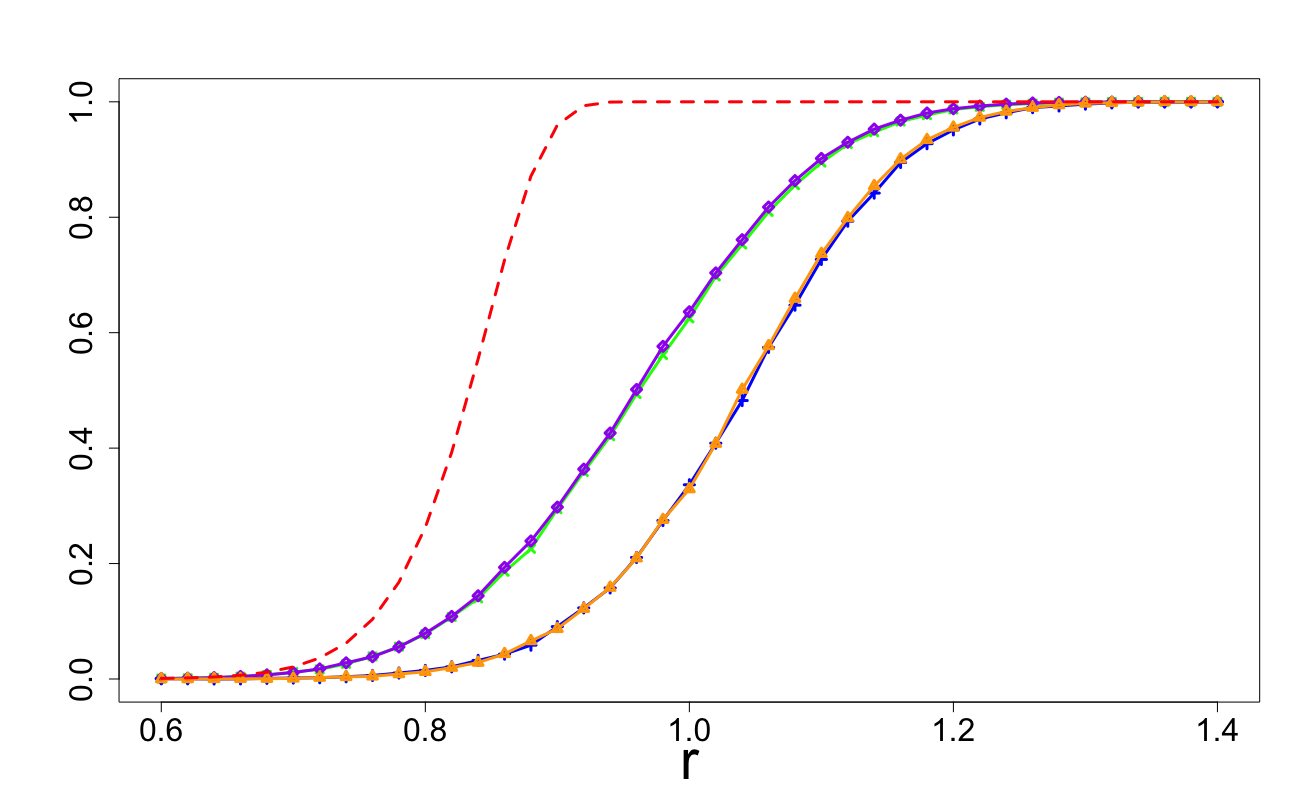}
  \caption{Sobol and $\delta$-cube points versus \\ the asymptotic covering:  $d=20, n=1024$. }
    \label{key_figure20}
\end{minipage}
\end{figure}

\AZ{  Figures~\ref{key_figure17}-\ref{key_figure18} are
 extended versions of  Fig.~\ref{key_figure1}.  Here we plot $F_d(r,X_n)$ as a function of $d$, where the radius is fixed from \eqref{normalised_radius} with $\gamma=0.1$ (the line $1-\gamma$ is depicted by a red solid line). For $X_n$ chosen uniformly in the cube $[0,1]^d$, we depict $F_d(r,X_n)$ with blue plusses. For $X_n$ chosen from a Sobol sequence in the whole cube $[0,1]^d$, we use  orange triangles. When the points in $X_n$ are uniform i.i.d.  within the $\delta$-cube with optimal $\delta$  we use green crosses. Finally, when points in $X_n$ are chosen from a Sobol sequence within the same $\delta$-cube we use purple diamond. Figures~\ref{key_figure17} and \ref{key_figure18}  illustrate two new key messages along with the message discussed in Figure~\ref{key_figure1}. Firstly, the use of low-discrepancy sequences seem to produce slightly better  results in comparison to random choice of points  for small dimensions but  in higher  dimensions the use of low-discrepancy  sequences (in our case, Sobol sequences) produces results that are almost equivalent to random sampling uniformly either in $[0,1]^d$ or in the optimally chosen $\delta$-cube. Secondly, in large dimensions sampling from a suitable  $\delta$-cube greatly outperforms the other schemes considered here being still far from the asymptotic results. These messages are further supported in Figures~\ref{key_figure19} and \ref{key_figure20}. Here we plot the asymptotic approximation $F_n(r)$ from Lemma~\ref{lem:1} (dashed red) and $F_d(r,X_n)$ as a function of $r$ for the following choices of  $X_n$: random in the cube $[0,1]^d$ (blue line with plusses),  chosen from a Sobol sequence in $[0,1]^d$ (orange line with triangles),  random in the $\delta$-cube with optimal delta (green line with crosses),   chosen from a Sobol sequence within the same $\delta$-cube (purple line with diamonds). We see that in Figure~\ref{key_figure19} for $d=10$, the Sobol sequence is slightly advantageous to the random uniform on the whole cube and $\delta$-cube for most interesting values of $\gamma$. Choosing $X_n$ as uniform within the $\delta$-cube produces better coverings than with Sobol's points in $[0,1]^d$ for most values of $\gamma$, but slightly worse for small $\gamma$. This slight advantage of the Sobol sequence in $[0,1]^d$ and in the $\delta$-cube diminishes in the case  $d=20$ shown in Figure~\ref{key_figure20}.
}

\begin{figure}[!h]
\centering
\begin{minipage}{.5\textwidth}
  \includegraphics[width=1\linewidth]{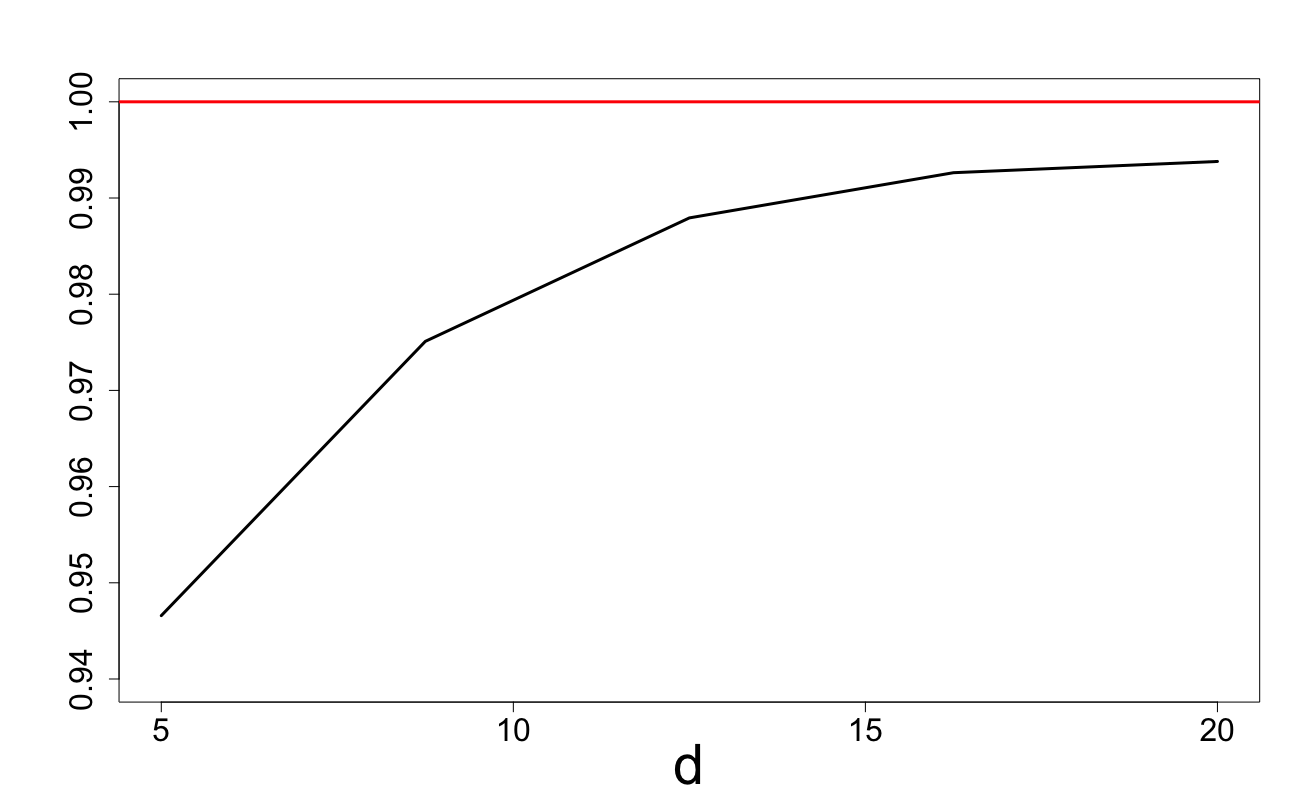}
  \caption{Efficiency of Sobol's points, $ n=2^{10}$. }
\label{sobol_efficiency1}
\end{minipage}%
\begin{minipage}{.5\textwidth}
  \centering
  \includegraphics[width=1\linewidth]{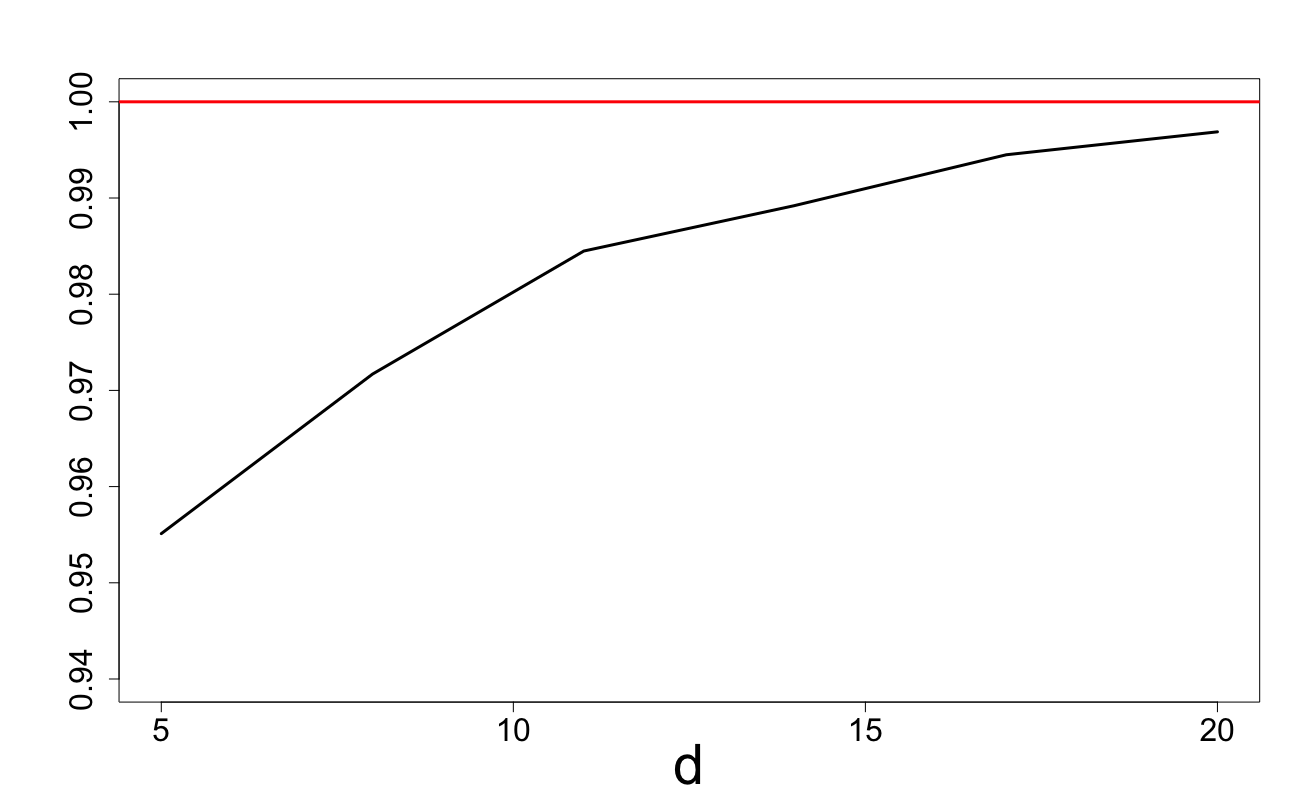}
  \caption{Efficiency of Sobol's points, $ n=2^{13}$. }
    \label{sobol_efficiency2}
\end{minipage}
\end{figure}

\AZ{
To further study the similarities in performance between $X_n$ chosen randomly in the $\delta$-cube with optimal $\delta$ and $X_n$ chosen from a Sobol sequence within the same $\delta$-cube, in Figures~\ref{sobol_efficiency1}-\ref{sobol_efficiency2} we plot the ratio of  the c.d.f.'s $F_d(r,X_{n,U})/F_d(r,X_{n,S})$ across different $d$. The subscript $U$ and $S$ respectively  differentiate between $X_n$ chosen randomly in the $\delta$-cube with optimal $\delta$ and $X_n$ chosen from a Sobol sequence within the same $\delta$-cube. For each value of $d$, $r$ is chosen so that $\max_{0\leq\delta\leq1}F_d(r,X_{n,U})=0.9$. 
}

\subsection{ \AZ{Non-uniform prior distribution for the target}}
\label{sec:delta1_3}

\AZ{In this section, we explore the effect a non-uniform prior distribution for $x_* \in \X$ has on the conclusions above formulated for the case of uniform distribution. We will assume that each component of $x_*$ has independent components  distributed according to the following symmetric beta distribution with density:
\bea
p_{\alpha}(t)= \frac{t^{\alpha-1}[1-t]^{\alpha-1}}{\mbox{Beta$(\alpha,\alpha)$}} \, , \text{ for some }\alpha> 0\,.
\eea
If $\alpha=1$, the density $p_{\alpha}(t)$ is uniform on $[0,1]$ while for $0<\alpha<1$ this density is U-shaped. In most cases below we choose the arcsine  density $p_{0.5}(t)$.
}

\AZ{
We then select the distribution of random points $x_j \in X_n$ to have similar shape  but constrained to the $\delta$-cube. More precisely, we assume that $x_j$ have independent components with the density
\bea
p_{\alpha,\delta}(t)= \frac{2\cdot (2\delta)^{1-2\alpha}}{\mbox{Beta$(\alpha,\alpha)$}} [\delta^2-(2t-1)^2]^{\alpha-1}\, , \;\;\frac{1-\delta}{2}<t<\frac{\delta+1}{2}\, ,\text{ for some }\alpha> 0\;\;{\rm and }\; 0\leq \delta\leq 1.
\eea
In the case $\alpha=1$, the points $x_j$ have uniform distribution on the cube $C_\delta$.
}

\begin{figure}[!h]
\centering
\begin{minipage}{.5\textwidth}
  \includegraphics[width=1\linewidth]{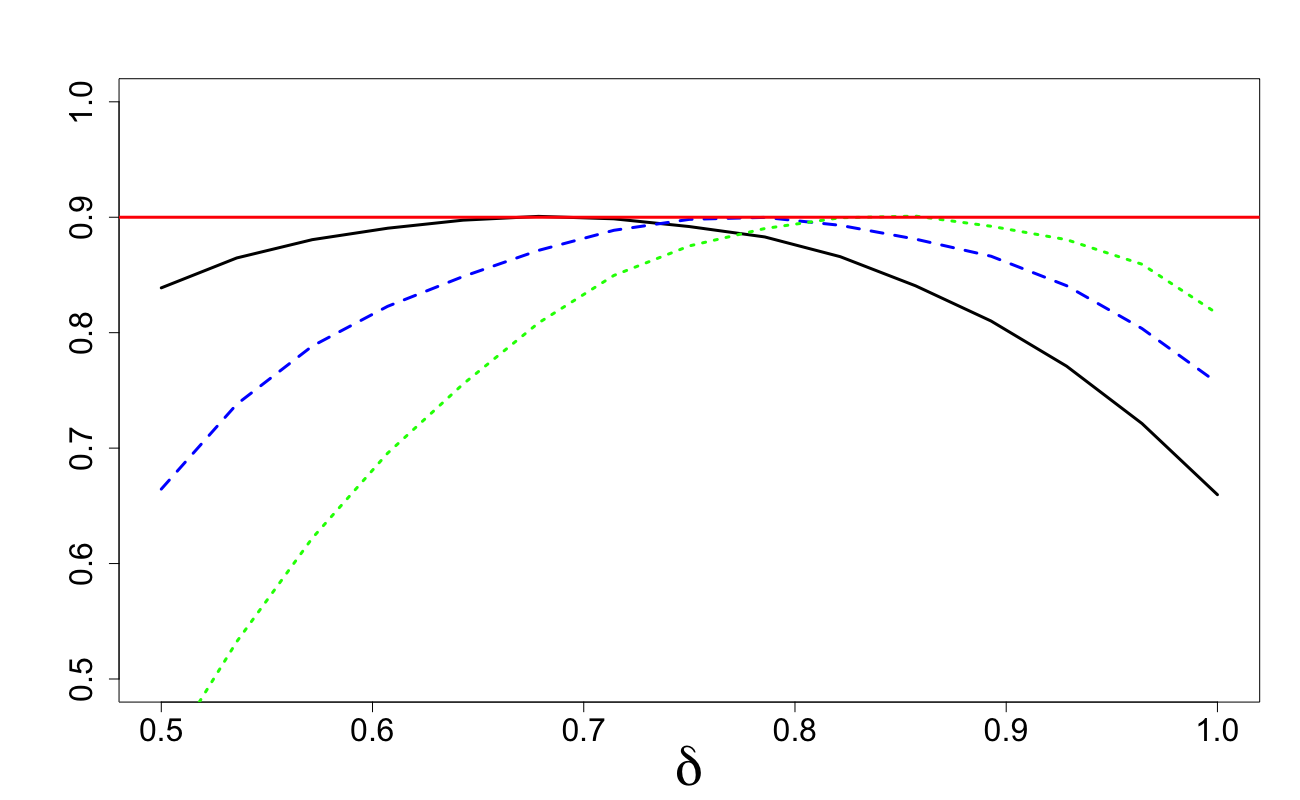}
  \caption{$d=20$ and $\alpha=0.5$ with \\ $n=1000,10000,100000$.  }
  \label{beta_sample}
\end{minipage}%
\begin{minipage}{.5\textwidth}
  \centering
  \includegraphics[width=1\linewidth]{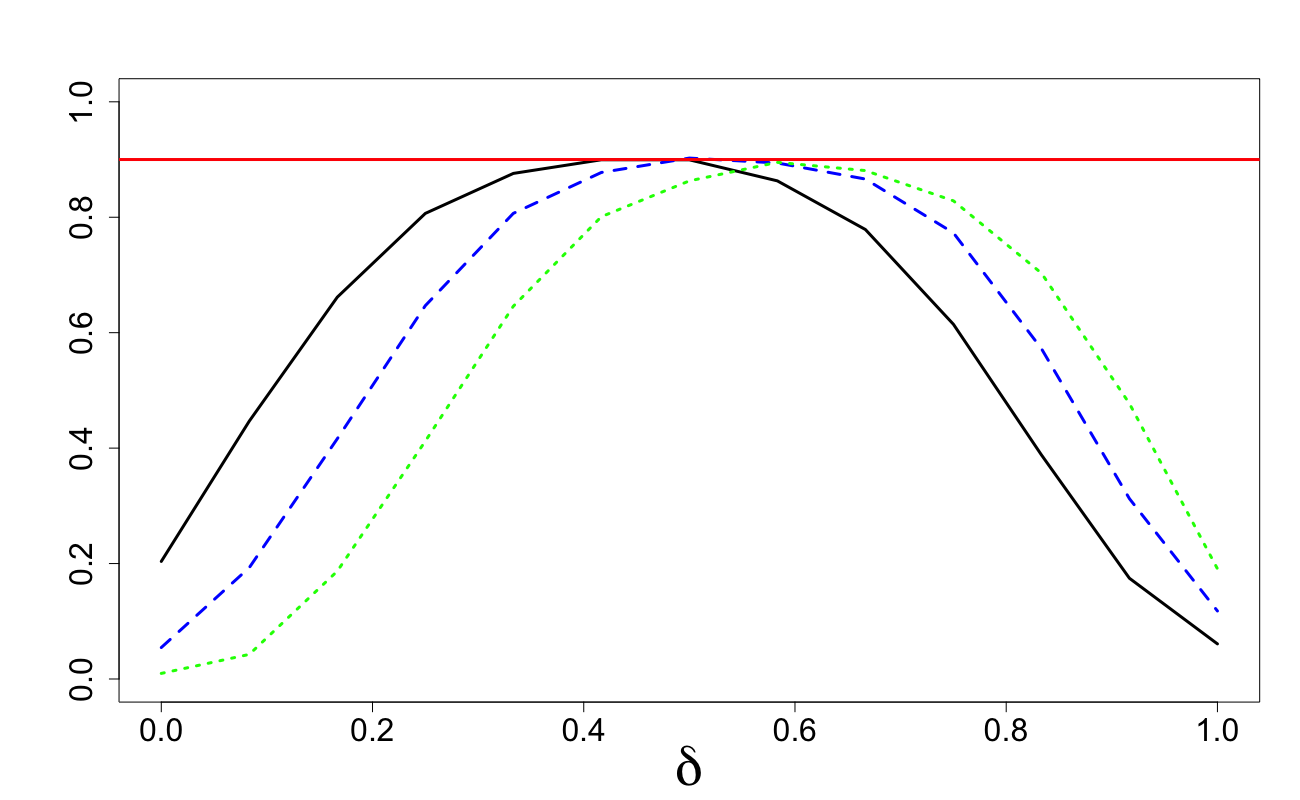}
  \caption{$d=50$ and $\alpha=0.5$ with  \\ $n=1000,10000,100000$. }
      \label{beta_sample2}
\end{minipage}
\end{figure}

\begin{figure}[!h]
\centering
\begin{minipage}{.5\textwidth}
  \includegraphics[width=1\linewidth]{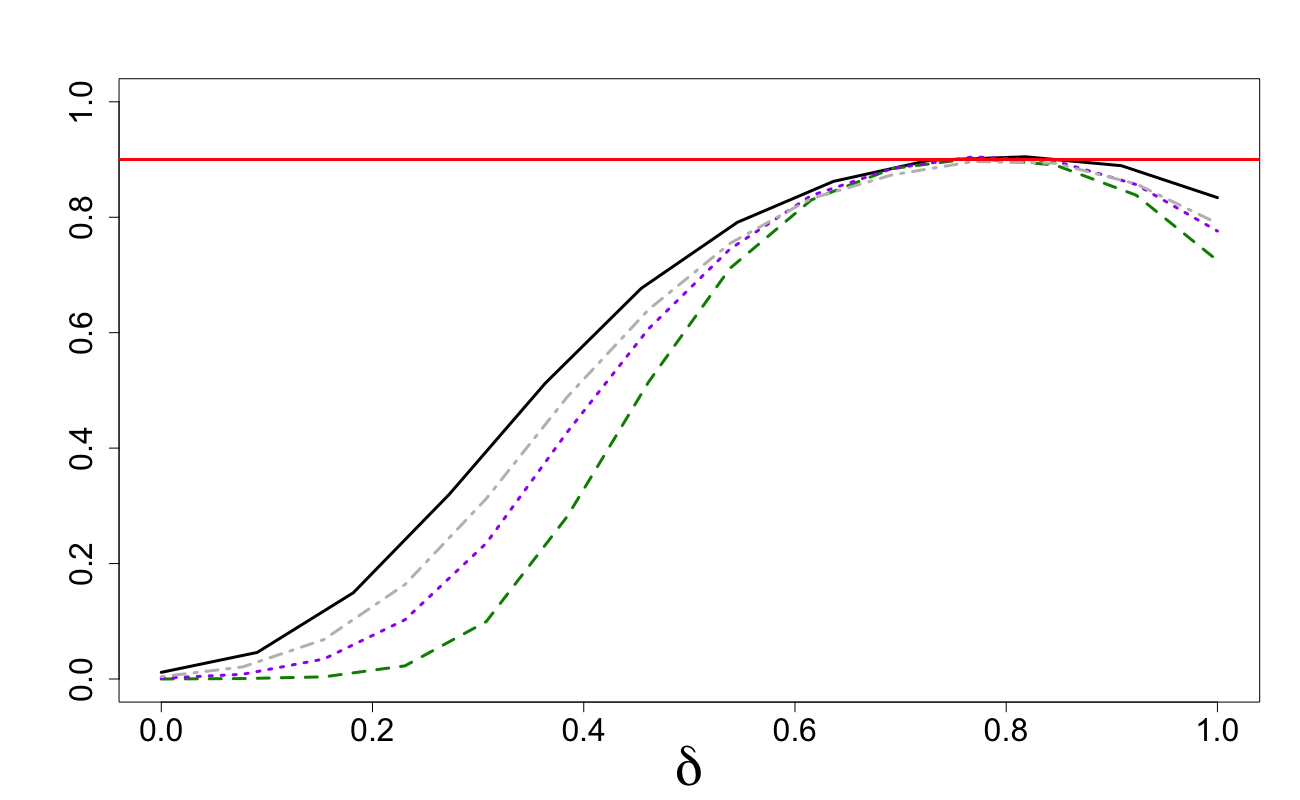}
  \caption{ $d=20, n=10000$ and \\$\alpha=0.25,0.5,0.75,1$. }
  \label{alpha1}
\end{minipage}%
\begin{minipage}{.5\textwidth}
  \centering
  \includegraphics[width=1\linewidth]{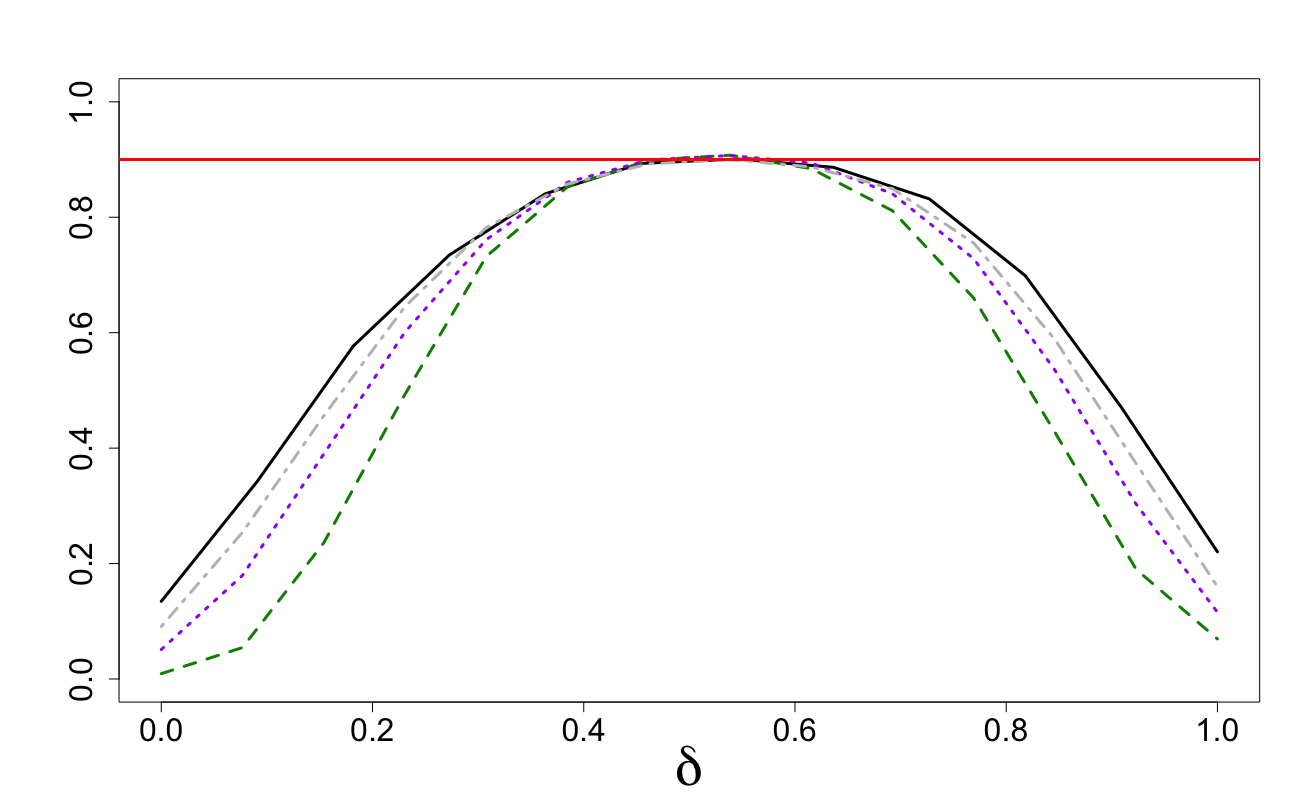}
  \caption{ $d=50, n=10000$ and\\ $\alpha=0.25,0.5,0.75,1$. }
      \label{alpha2}
\end{minipage}
\end{figure}

\AZ{
 Figures~\ref{beta_sample}--\ref{beta_sample2} are similar to  Figures~\ref{key_figure13}--\ref{key_figure15}, but with the key difference of assuming a non-uniform prior distribution for $x_*$. For different values of $d$ and $ n$ and for $\alpha=0.5$, we plot $F_d(r,X_n)$ as a function of $\delta$. For each $d$ and $n$, the value of $r$ has been chosen so that $\max_{0\leq\delta\leq1}F_d(r,X_n)=0.9$. In these figures, the values of $F_d(r,X_n)$ for $n=1000, 10000, 100000$ are shown with a solid black line, dashed blue line and dotted green line respectively.   Figures~\ref{alpha1}-\ref{alpha2} are  similar  to Figures~\ref{beta_sample}--\ref{beta_sample2}, but with varying  values of $\alpha$ and fixed $n=10000$. In these figures, we selected $\alpha=0.25$ (dashed dark green), $\alpha=0.5$ (dotted purple), $\alpha=0.75$ (dot-dashed grey) and $\alpha=1$ which gives the uniform distribution  (solid black). Figures~\ref{beta_sample}--\ref{beta_sample2} clearly demonstrate that the `$\delta$-effect' is still significant in this non-uniform setting.  
}

\section{Intersection of one ball with the cube}\label{Petrov_section}

As a result of \eqref{eq:prod}, our main quantity of interest in this section will be the probability
\be \label{eq:prob_u}
P_{U,\delta,r}:=\mathbb{P}_{_X} \left\{ \| U\!-\!X \|\! \leq \! { r } \right\}\!= \! \mathbb{P}_{_X} \left\{ \| U\!-\!X \|^2 \leq  { r^2 } \right\}\!= \! \mathbb{P} \left\{\sum_{j=1}^d (u_j\!-\!x_j)^2 \leq  { r }^2 \right\}   \;\;
\ee
in the case when $X$ has  the uniform distribution on the $\delta$-cube $[1/2-\delta/2, 1/2+\delta/2]^d$  and ${U=(u_1, \ldots, u_d) \in \mathbb{R}^d}$ is fixed. The case of $\delta=1$ will be directly applicable to Section~\ref{sec:jensen}. Because of the results of Section~\ref{sec:jensen}, we will bear in mind two typical choices of $U$ will be $U=\bm{1/2}$ and $U=\bm{3/4}$ but will formulate results for general $U$.

For fixed $u\in \mathbb{R}^d$, consider the r.v.  $\eta_{u,\delta} = (z-u)^2$, where $z$ has density
\be\label{Beta_density}
p_{\delta}(t)= 1/\delta\, , \;\;(1-\delta)/2<t<(1+\delta)/2\, ,\text{ for some } 0\leq \delta\leq 1.
\ee
The first three  central moments of $\eta_{u,\delta}$  are:
\be
\label{1eq:inters1c}
\mu_{u}^{(1)} &=& \mathbb{E}\eta_{u,\delta} =\left(u-\frac 12 \right )^2+ \frac {{{\delta}}^{2}}{12} \, ,\\
\label{1eq:inters1c2}
\mu_{u}^{(2)}&=&{\rm var} (\eta_{u,\delta}) =
{\frac {\delta^{2}}{3}}
 \left[\left(u-\frac 12 \right )^2+
 {\frac {{{\delta}}^{2}}{60}}
  \right]  \, ,\\
\label{1eq:inters1c3}
\mu_{u}^{(3)}&=&\mathbb{E} \left[\eta_{u,\delta} - \mu_{u}^{(1)}\right]^3 =
{\frac {4 {{\delta}}^{4}}{ 15 }} \left[  \left(u-\frac 12 \right )^2+ {\frac {{{\delta}}^{2} }{252}}
 \right]
\, .
\ee

Then for given $U=(u_1, \ldots, u_d) \in \mathbb{R}^d$, consider the random variable
\bea
\| U-X \|^2 =\sum_{i=1}^d \eta_{u_i,\delta,\alpha} =\sum_{j=1}^d (u_j-x_j)^2\, ,
\eea
where we assume that $X=(x_1, \ldots, x_d) $ is a random vector with i.i.d. components $x_i$ with density \eqref{Beta_density}.
From \eqref{1eq:inters1c}, its mean is
\bea
\mu=\mu_{d,\delta,U}:=\mathbb{E}\| U-X \|^2  =\|U-\bm{1/2}\|^2 +\frac {{d{\delta}}^{2}}{12}\, .
\eea
Using independence of $x_1, \ldots, x_d$ and  \eqref{1eq:inters1c2}, we obtain
\bea
 {\sigma}_{d,\delta,U}^2 :={\rm var}(\| U-X \|^2 )  =
{\frac {\delta^{2}}{3}}
 \left[\|U- \bm{1/2}\|^2+
 {\frac {{d{\delta}}^{2}}{ 60}}
  \right]  \, ,
\eea
and from independence of $x_1, \ldots, x_d$ and  \eqref{1eq:inters1c3} we get
\be\label{third_central}
  {\mu}_{d,\delta,U}^{(3)} := \mathbb{E}\left[\| U-X \|^2- \mu\right]^3  = \sum_{j=1}^d     \mu_{u_j}^{(3)} = {\frac {\,{{\delta}}^{4}}{ 15 }} \left[  \|U-\bm{1/2}\|^2+ {\frac {{d{\delta}}^{2} }{
252 }}
 \right]
\, .\;\;\;\;\;\;
\ee

If $d$ is large enough then the conditions of the CLT for $\| U-X \|^2$ are approximately met and  the distribution of $\| U-X \|^2 $
 is approximately normal with mean $\mu_{d,\delta,U}$ and variance ${\sigma}_{d,\delta,U}^2$. That is, we can approximate the probability
$P_{U,\delta,r}= \mathbb{P}_{_X} \left\{ \| U\!-\!X \|\! \leq \! { r } \right\}$
by
\be
\label{1eq:inters2f1}
P_{U,\delta,r}\!\cong \Phi \left(\frac{{ r }^2-\mu_{d,\delta,U}}{{\sigma}_{d,\delta,U}} \right) \, ,
\ee
where $\Phi (\cdot)$   is the c.d.f. of the standard normal distribution:
$$
\Phi (t) = \int_{-\infty}^t \varphi(v)dv\;\;{\rm with}\;\; \varphi(v)=\frac{1}{\sqrt{2\pi}} e^{-v^2/2}\, .
$$
The approximation \eqref{1eq:inters2f1} has acceptable accuracy if the probability $P_{U,\delta,r}$ is not very small; for example, it falls inside a $2\sigma$-confidence interval generated by the standard normal distribution.

To improve on the usual CLT approximation, we use Edgeworth-type  expansion in the CLT for sums of independent  {non-identically} distributed r.v.
by
 V.Petrov, see \cite{petrov2012sums}:
\be\label{petrov_approx}
P\left(\frac{\| U-X \|^2-\mu_{d,\delta,U}}{\sigma_{d,\delta,U}} \leq t \right) = \Phi(t)
+\sum_{\nu=1}^{\infty}\frac{Q_{\nu,d}(t)}{d^{\nu/2}} \, ,
\ee
where
\bea
Q_{\nu,d}(t) = -\varphi(t)\sum H_{\nu+2s-1}(t)\prod_{m=1}^{\nu}\frac{1}{k_m!}\left(\frac{\lambda_{m+2,d}}{(m+2)!} \right)^{k_m} \,,
\eea
\bea
\lambda{}_{\nu,d} = \frac{d^{(\nu-2)/2}}{\sigma_{d,\delta,U}^\nu}\sum_{j=1}^{d}\gamma_{\nu,j} \,
\eea
$\gamma_{\nu,j}$ is the  cumulant of order $\nu$ at $(u_j-x_j)^2$, $H_m$ is the Chebyshev-Hermite polynomial of degree $m$ and the summation is carried out over all non-negative integer solutions of the equation
\bea
k_1+2k_2+\cdots+\nu k_\nu =\nu \\
s=k_1+k_2+\cdots+k_\nu \,.
\eea
The partition function $p(\nu)$ provides the number of possible partitions of a non-negative integer $\nu$ and therefore at each value of $\nu$ provides the number of terms in the summation. The sequence has the generating function
\bea
\sum_{\nu=0}^{\infty} p(\nu)x^{\nu} = \prod_{k=1}^{\infty} \left(\frac{1}{1-x^k} \right) ;
\eea
the first few values are:
1, 1, 2, 3, 5, 7, 11, 15, 22, 30, 42, 56, 77, 101. The first few terms in the summation (including Hermite polynomials) are provided in \cite[p. 139]{Petrov}.

In the case of $U=\bm{1/2}$ or $U=\bm{3/4}$, the random variables $(u_j-x_j)^2$ will be i.i.d. For this case $\lambda_{\nu,d}$ does not depend on $d$ and thus we have the slight simplification $\lambda{}_{\nu,d}=\lambda{}_{\nu} = \gamma_\nu/\sigma_{d,\delta,U}^\nu$.

\begin{figure}[!h]
\centering
\begin{minipage}{.5\textwidth}
  \includegraphics[width=1\linewidth]{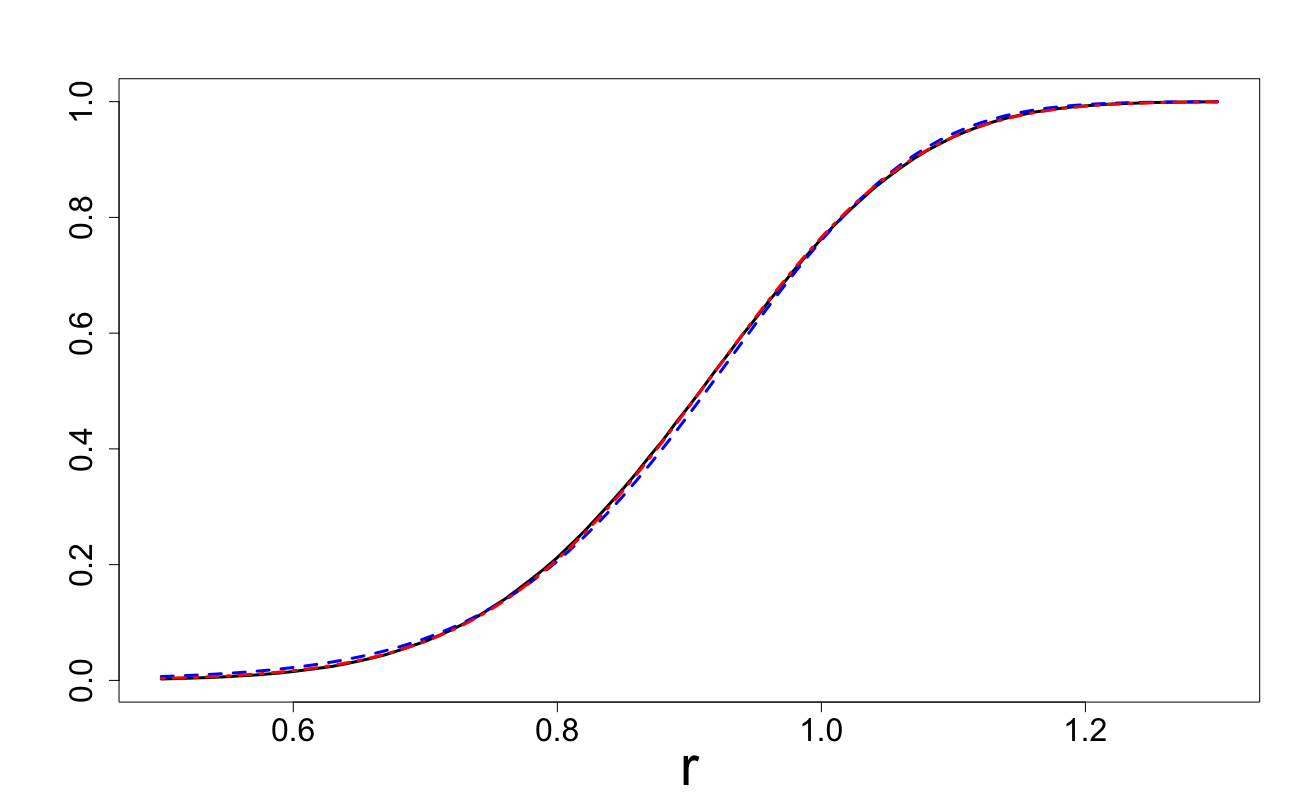}
  \caption{$d=10, U=\bm{1/2}$. }
  \label{key_figure22}
\end{minipage}%
\begin{minipage}{.5\textwidth}
  \centering
  \includegraphics[width=1\linewidth]{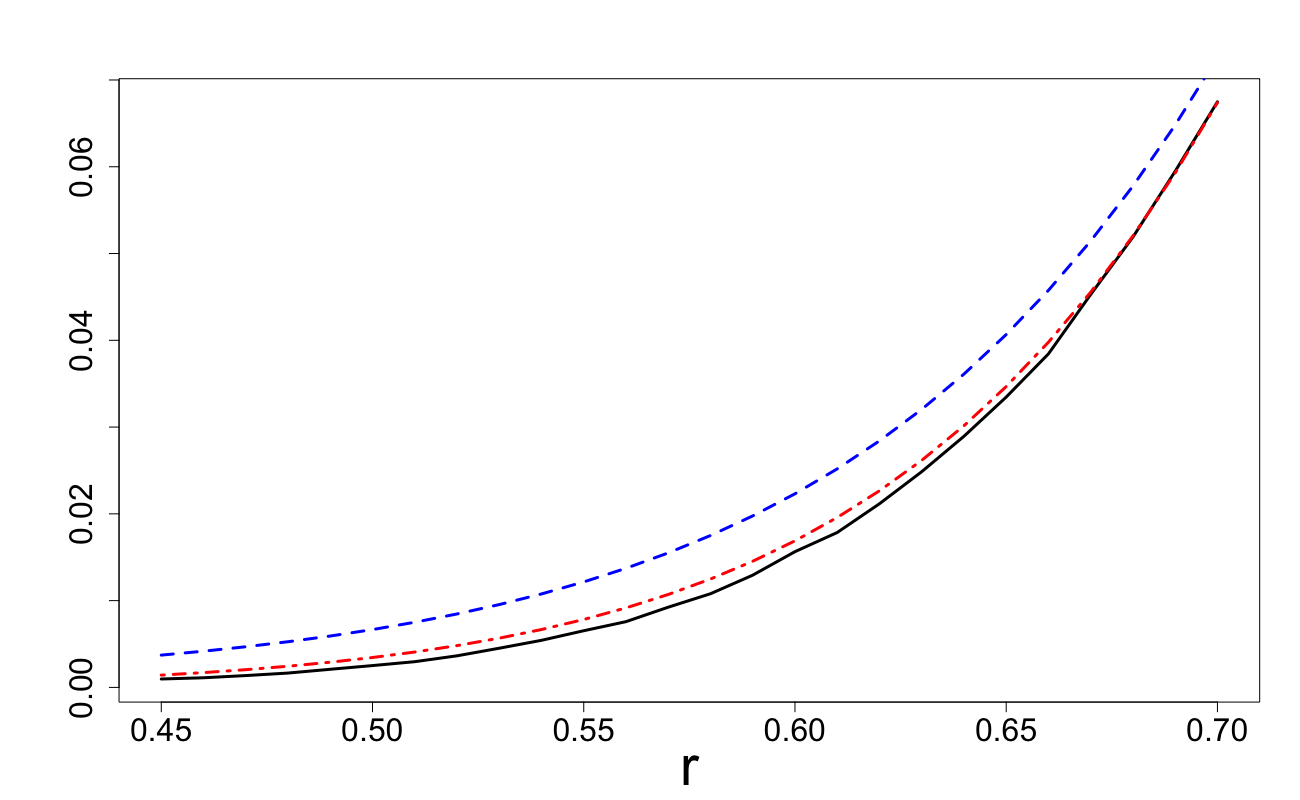}
  \caption{$d=10, U=\bm{1/2}$. }
    \label{key_figure23}
\end{minipage}
\end{figure}

\begin{figure}[!h]
\centering
\begin{minipage}{.5\textwidth}
  \includegraphics[width=1\linewidth]{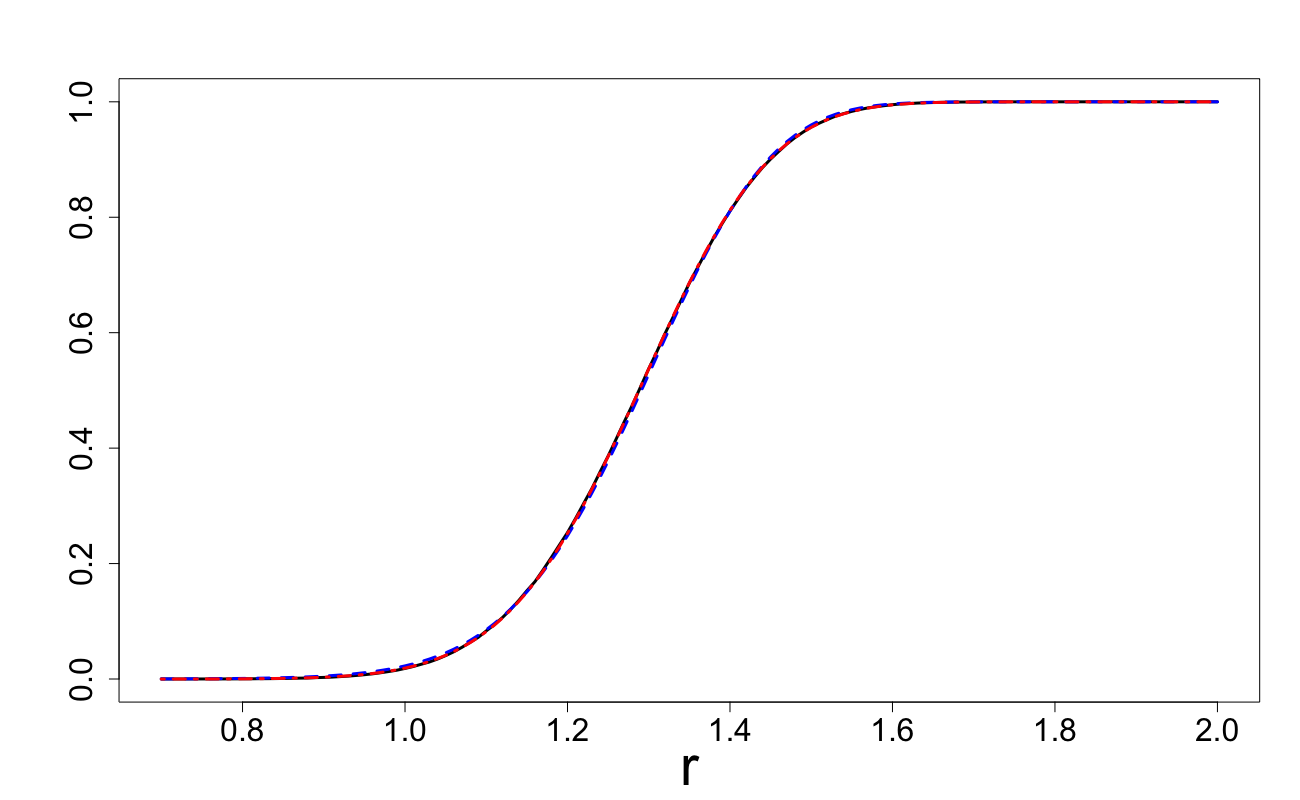}
  \caption{$d=20, U=\bm{1/2}$. }
  \label{key_figure222}
\end{minipage}%
\begin{minipage}{.5\textwidth}
  \centering
  \includegraphics[width=1\linewidth]{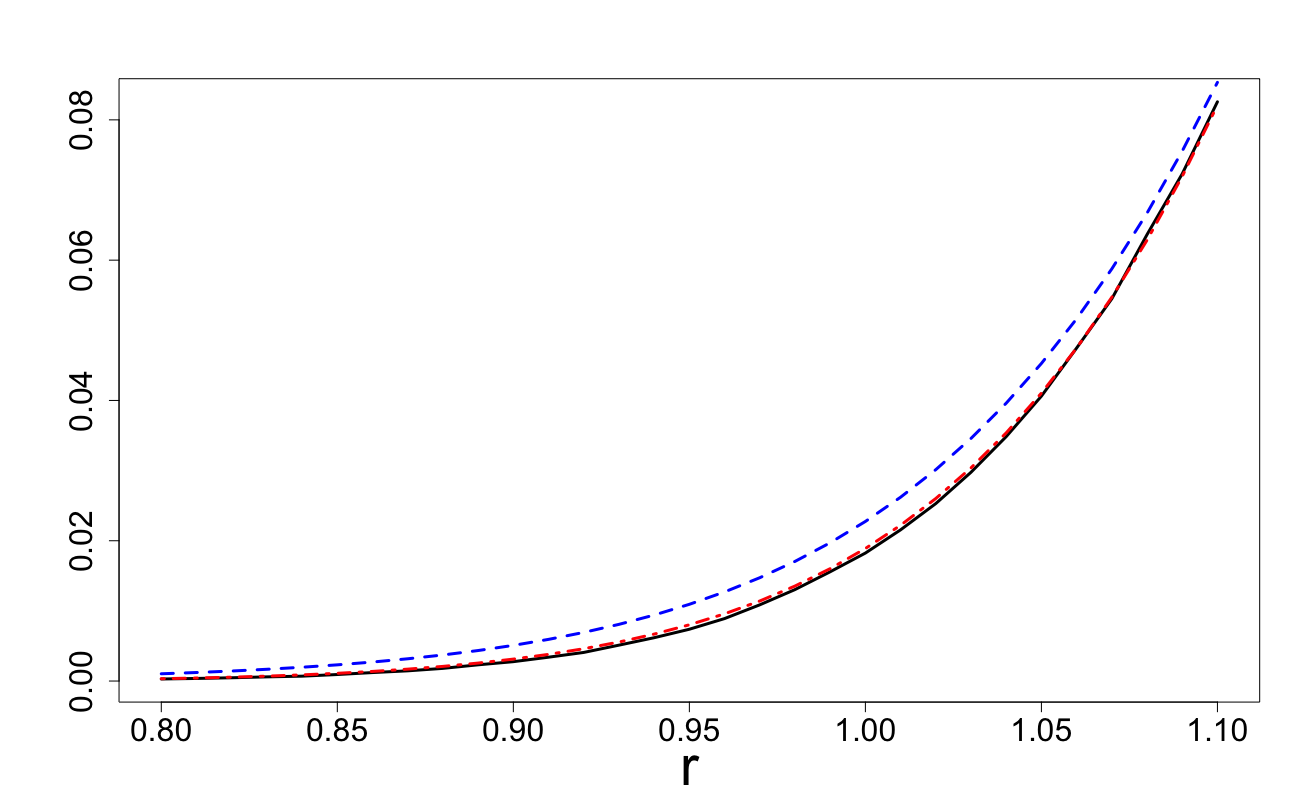}
  \caption{$d=20, U=\bm{1/2}$. }
    \label{key_figure232}
\end{minipage}
\end{figure}

\begin{figure}[!h]
\centering
\begin{minipage}{.5\textwidth}
  \includegraphics[width=1\linewidth]{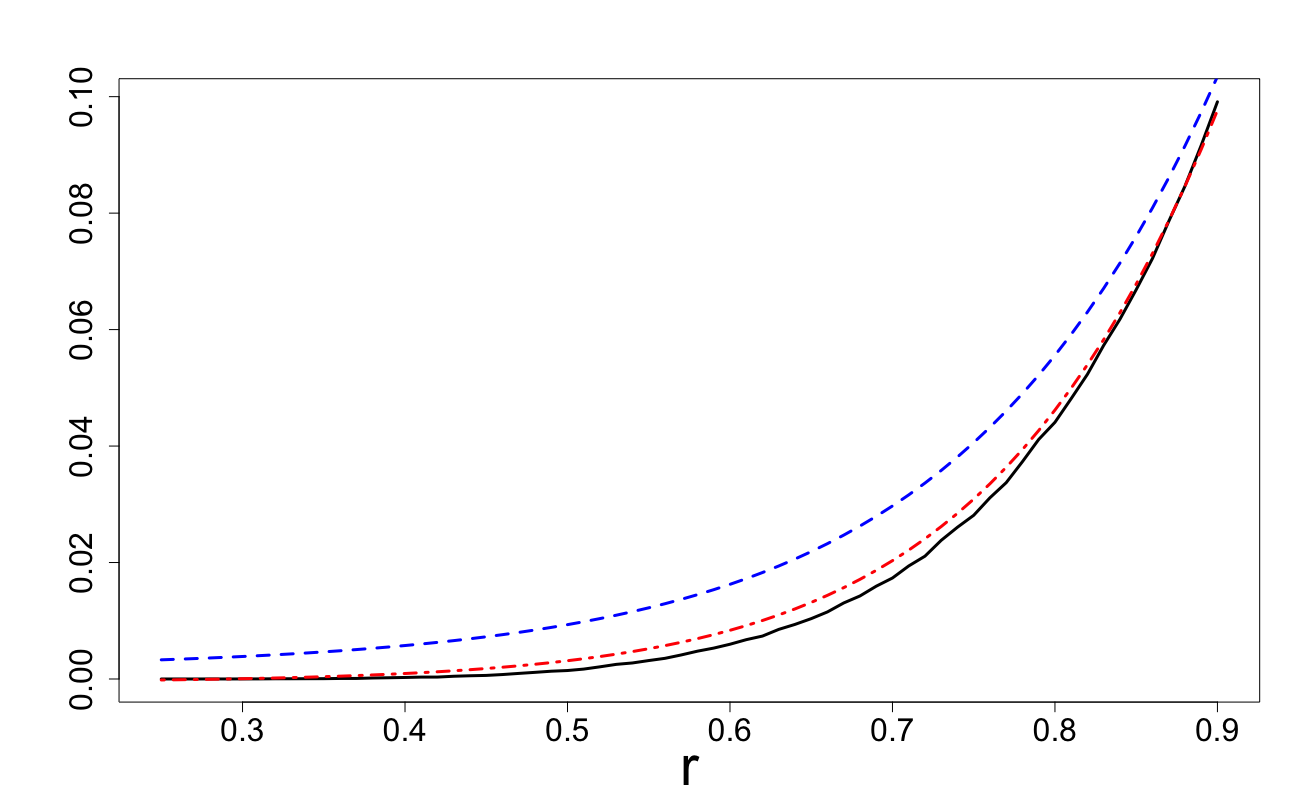}
  \caption{$d=10, U=\bm{3/4}$. }
  \label{key_figure24}
\end{minipage}%
\begin{minipage}{.5\textwidth}
  \centering
  \includegraphics[width=1\linewidth]{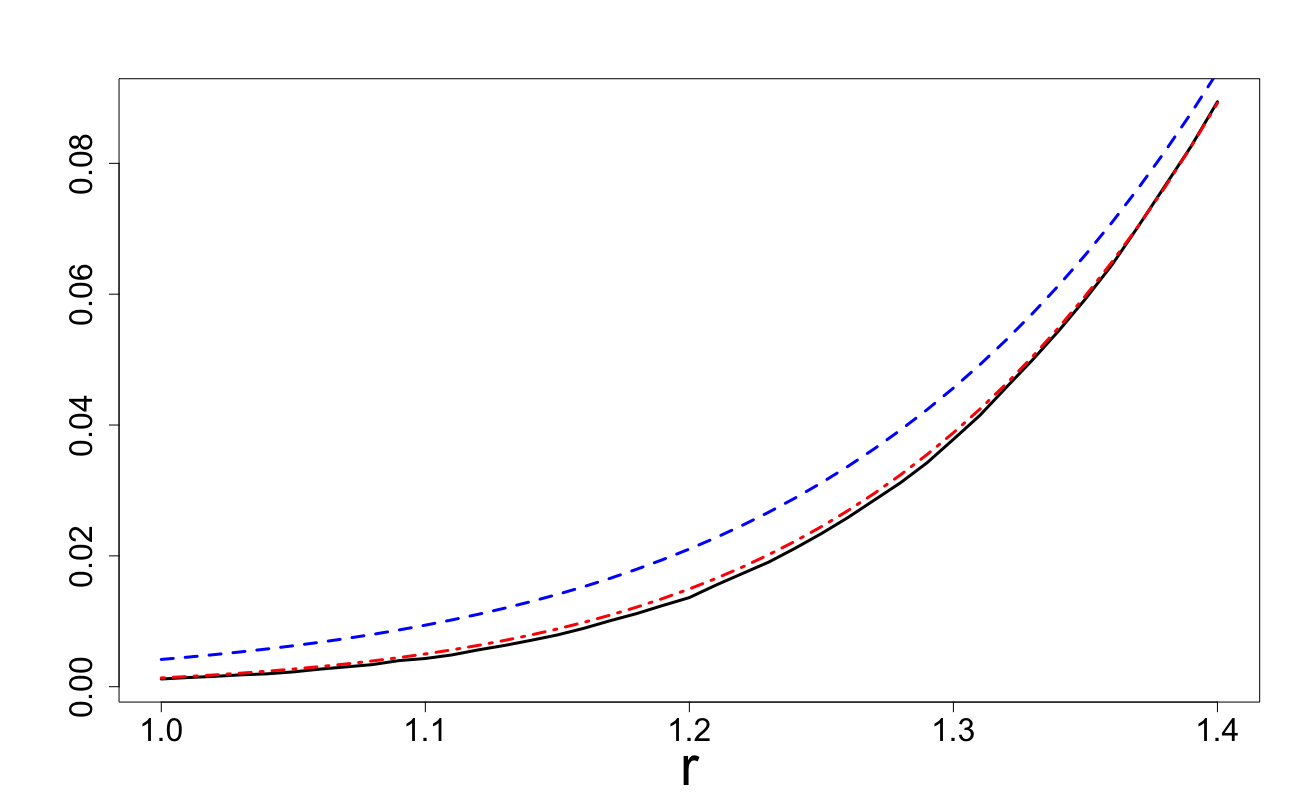}
  \caption{$d=20, U=\bm{3/4}$. }
    \label{key_figure25}
\end{minipage}
\end{figure}

In Figures~\ref{key_figure22}--\ref{key_figure25}, we plot $P_{U,1,r}$ for $U=\bm{1/2}$ and $U=\bm{3/4}$ as a function of $r$ with a solid black line. In these figures, we demonstrate the accuracy of approximation \eqref{1eq:inters2f1} with a dashed blue line. With a dot-dashed red line, we plot the accuracy of an approximation obtained by taking one additional term in the expansion given in \eqref{petrov_approx}; this requires use of the third central moment given in \eqref{third_central}. We can see that overall, for $d=10$ and $d=20$, the approximations are fairly accurate. However, when considering covering by $n$ balls it is more important to focus on the lower tail. Figures~\ref{key_figure23}, \ref{key_figure232}, \ref{key_figure24} and \ref{key_figure25} demonstrate that taking one additional term in the Petrov's expansion \eqref{petrov_approx} produces a significant improvement in accuracy.

\section{Conclusions}

\label{sec:conclusions}

We have  considered continuous global optimization problems, where the feasible region $\X$ is a compact subset of $\mathbb{R}^d$.
As a strategy for exploration, we have mostly considered  sampling of i.i.d. random points either in $\X$ or a suitable subset of $\X$.

We have distinguished between between `small', `medium' and `high' dimensional problems depending on the following relations between $d$ and $n_{\max}$ (which is the maximum possible number of points available for space exploration):
\begin{itemize}
  \item[(S)] small dimensions:  $n_{\max} \gg 2^d$ (roughly, $d < 10$);
  \item[(M)] medium dimensions: $n_{\max} $ is comparable to $ 2^d$   (roughly, $10 \leq d \leq 20$);
  \item[(H)] high dimensions: $n_{\max} \ll 2^d$ (roughly, $d > 20$).
\end{itemize}

We only considered the situations (M) and (H), where we have demonstrated the following effects:  (i)  the actual convergence of randomized exploration schemes is much slower than that given by the classical estimates, which are based on the asymptotic properties of random points; (ii)  the usually recommended space exploration schemes  are practically inefficient as the asymptotic regime is unreachable. In particular, we have shown: (ii-a) uniform sampling on entire~$\X$ is much less efficient than uniform sampling on a suitable subset of $\X$,  and (ii-b)  the effect of replacement of  random points by a low-discrepancy sequence  is very small so that using low-discrepancy sequences and other deterministic constructions does not lead to significant improvements
(unless the number of evaluation points $n=n_{\max}$ is fixed to some particular value like $2^{d} $ or $2^{d-1}$, see \cite{noonan2022efficient}). We believe that the effects (i) and  (ii) have not been stated in literature, at least in this generality. The effect (ii-a) has been numerically demonstrated in our previous papers \cite{us,second_paper}. The effect (ii-b) enhances one of the main messages of the paper   \cite{pepelyshev2018performance}.

It was not the purpose of the paper to give the most effective exploration schemes. However, the results of this paper, along with  studies reported in \cite{us,second_paper} and \cite{noonan2022efficient},
allow us to give several general recommendations on efficient organization of exploration strategies in the situations (M) and (H), at least when $\X$ is a cube.

In a high-dimensional cube $\X=[0,1]^d$ with $d>20$ and $1<n_{\max}<2^{d}$, we propose the following strategy of construction of nested exploration designs $X_n$: $x_1=\bm{1/2}$ (the centre of  $\X$) and the other points $x_j$ are taken randomly among the vertices of a cube $[1/4,3/4]^d$. Sampling from vertices can be done without replacement (see \cite{second_paper}) and, moreover, we can keep the points $x_j$ so that the Hamming distance between them is at least
  {  $\lfloor d-\log_2 (n_{\max}-1) \rfloor+1$. }

It is more difficult to be so specific in the situation (M) as there may be different relations between $2^d$, $n_{\min}$ and $n_{\max}$. A good strategy would be using the product of arcsine distributions on a suitable $\delta$-cube  (see \cite{second_paper}); this distribution is slightly superior to the uniform on $\delta$-cube
of Section~\ref{sec:delta1} (with different values of $\delta$ optimized for the respective distribution). An even more natural strategy would be sampling in $2^d$ small cubes (or side-length $\varepsilon$) surrounding the vertices of a $\delta$-cube   $C_\delta=[{ 1/2-{\delta}/2, 1/2+{\delta}/2}]^d $ (after placing $x_1=\bm{1/2}$).
The choice of $\delta$ and $\varepsilon$ depends on $2^d$, $n_{\min}$ and $n_{\max}$ are requires a separate study.
As usual, reduction of randomness in sampling makes any of these schemes marginally more efficient.

\section*{Data availability statement}
Data sharing not applicable to this article as no datasets were generated or analysed during the current study.

\bibliographystyle{plainnat}

\bibliography{large_dimension_enc_1}

\end{document}